# Boundary rigidity and stability for generic simple metrics

Plamen Stefanov*
Department of Mathematics
Purdue University
West Lafayette, IN 47907, USA

Gunther Uhlmann[†]
Department of Mathematics
University of Washington
Seattle, WA 98195, USA

**Abstract**

We study the boundary rigidity problem for compact Riemannian manifolds with boundary $(M, g)$: is the Riemannian metric $g$ uniquely determined, up to an action of diffeomorphism fixing the boundary, by the distance function $\rho_g(x, y)$ known for all boundary points $x$ and $y$? We prove in this paper global uniqueness and stability for the boundary rigidity problem for generic simple metrics. More specifically, we show that there exists a generic set $\mathcal{G}$ of simple Riemannian metrics and an open dense set $\mathcal{U} \subset \mathcal{G} \times \mathcal{G}$, such that any two Riemannian metrics in $\mathcal{U}$ having the same distance function, must be isometric. We also prove Hölder type stability estimates for this problem for metrics which are close to a given one in $\mathcal{G}$.

## 1  Introduction

Let $(M, g)$ be a Riemannian manifold with boundary. Denote by $\rho_g$ the distance function in the metric $g$. We consider the inverse problem of whether $\rho_g(x, y)$, known for all $x$, $y$ on $\partial M$, determines the metric uniquely. This problem arose in geophysics in an attempt to determine the inner structure of the Earth by measuring the travel times of seismic waves. It goes back to Herglotz [H] and Wiechert and Zoeppritz [WZ]. Although the emphasis has been in the case that the medium is isotropic, the anisotropic case has been of interest in geophysics since it has been found that the inner core of the Earth exhibits anisotropic behavior [Cr]. In differential geometry this inverse problem has been studied because of rigidity questions and is known as the boundary rigidity problem. It is clear that one cannot determine the metric uniquely. Any isometry which is the identity at the boundary will give rise to the same measurements. Furthermore, the boundary distance function only takes into account the shortest paths and it is easy to find counterexamples to unique determination, so one needs to pose some restrictions on the metric. Michel [Mi], conjectured that a *simple* metric $g$ is uniquely determined, up to an action of a diffeomorphism fixing the boundary, by the boundary distance function $\rho_g(x, y)$ known for all $x$ and $y$ on $\partial M$. We recall

**Definition 1** *We say that the Riemannian metric $g$ is* simple *in $M$, if $\partial M$ is strictly convex w.r.t. $g$, and for any $x \in M$, the exponential map $\exp_x : \exp_x^{-1}(M) \to M$ is a diffeomorphism.*

Note that a simple metric $g$ in $M$ can be extended to a simple metric in some $M_1$ with $M \subset\subset M_1$. If we fix $x = x_0 \in M$ above, we also obtain that each simple manifold is diffeomorphic to a (strictly convex)

*Partly supported by NSF Grant DMS-0196440
[†]Partly supported by NSF and a John Simon Guggenheim fellowship

domain $\Omega \subset \mathbf{R}^n$ with the Euclidean coordinates $x$ in a neighborhood of $\Omega$ and a metric $g(x)$ there. For this reason, it is enough to prove our results for domains $\Omega$ in $\mathbf{R}^n$.

Unique recovery of $g$ (up to an action of a diffeomorphism) is known for simple metrics conformal to each other [C1], [B], [Mu1], [Mu2], [MuR], [BG], for flat metrics [Gr], for simple locally symmetric spaces of negative curvature [BCG]. In two dimensions in [C2] and [O] it was known for simple metrics with negative curvature recently it was shown in [PU] for simple metrics with no restrictions on the curvature. In [SU2], the authors proved this for metrics in a small neighborhood of the Euclidean one. This result was used in [LSU] to prove a semiglobal solvability result.

It is known [Sh1], that a linearization of the boundary rigidity problem near a simple metric $g$ is given by the following integral geometry problem: show that if for a symmetric tensor of order 2, which in coordinates is given by $f = (f_{ij})$, the geodesic X-ray transform

$$I_g f(\gamma) = \int f_{ij}(\gamma(t)) \dot{\gamma}^i(t) \dot{\gamma}^j(t) \, \mathrm{d}t$$

vanishes for all geodesics $\gamma$ in $M$, then $f = dv$ for some vector field $v$ with $v|_{\partial M} = 0$, where $dv$ denotes the symmetric differential

$$[dv]_{ij} = \frac{1}{2} \left( \nabla_i v_j + \nabla_j v_i \right), \tag{1.1}$$

and $\nabla_k v$ denote the covariant derivatives of the vector field $v$.

We will refer to this property as *s-injectivity* of $I_g$. More precisely, we have.

**Definition 2** *We say that $I_g$ is s-injective in $M$, if $I_g f = 0$ and $f \in L^2(M)$ imply $f = dv$ with some vector field $v \in H_0^1(M)$.*

If $f = dv$ above, then supp $v$ is included in the smallest simply connected closed set containing supp $f$, see Proposition 3.

On the other hand, it is easy to see that $I_g dv = 0$ for any such $v$. This is the linear version of the fact that the $\rho_g$ does not change on $(\partial M)^2 := \partial M \times \partial M$ under an action of a diffeomorphism as above. The inversion of $I_g$ is a problem of independent interest in integral geometry, and our first two theorems are related to it. S-injectivity of $I_g$ was proved in [PS] for metrics with negative curvature, in [Sh1] for metrics with small curvature and in [ShU] for Riemannian surfaces with no focal points. A conditional and non-sharp stability estimate for metrics with small curvature is also established in [Sh1]. This estimate was used in [CDS], [E] to get local uniqueness results for the boundary rigidity problem under the same condition. In [SU3], we proved stability estimates for s-injective metrics (see (1.2) and section 2 below) and sharp estimates about the recovery of a 1-form $f = f_j dx^j$ and a function $f$ from the associated $I_g f$. The stability estimates proven in [SU3], were used to prove local uniqueness for the boundary rigidity problem near any simple metric $g$ with s-injective $I_g$.

Before stating our results we give a definition, see also [Tre].

**Definition 3** *We say that $f$ is analytic in $K \subset \mathbf{R}^n$ (not necessarily open), and denote $f \in \mathcal{A}(K)$, if there exists an open set $U \supset K$, such that $f$ extends analytically in $U$.*

Let $\Omega$ be a bounded domain in $\mathbf{R}^n$, $n \geq 2$ with smooth boundary. We show in section 4.

**Theorem 1** *Let $g$ be a simple metric in $\Omega$, real analytic in $\bar{\Omega}$. Then $I_g$ is s-injective.*

The proof of Theorem 1 is based on the following. For analytic metrics, the normal operator $N_g = I_g^* I_g$, where $I_g^*$ denotes the operator adjoint to $I_g$ with respect to an appropriate measure (see section 2), is an analytic pseudodifferential operator with a non-trivial null space. In section 3 we construct an analytic parametrix that allows to reconstruct the solenoidal part of a tensor field from its geodesic X-ray transform. As shown in [SU3], this implies a stability estimate for $I_g$ for analytic simple $g$. In next theorem we show something more, namely that we have a stability estimate for $g$ in a neighborhood of each analytic metric, which leads to stability estimates for generic metrics.

Let $M_1 \supset M$ be a compact manifold which is a neighborhood of $M$ and $g$ extends as a simple metric there. We always assume that our tensors are extended as zero outside $M$, which may create jumps at $\partial M$. In (5.1), (5.2), see also [SU3], we define the space $\widetilde{H}^2(M_1)$ that in particular satisfies $H^1(M_1) \subset \widetilde{H}^2(M_1) \subset H^2(M_1)$. On the other hand, $f \in H^1(M)$ implies $N_g f \in \widetilde{H}^2(M_1)$ despite the possible jump of $f$ at $\partial M$. It is known that every symmetric tensor $f \in L^2$ can be uniquely decomposed into a *solenoidal* part $f^s$, and a *potential* part $dv$ as above, i.e., $f = f^s + dv$, see section 2.

**Theorem 2** *There exists $k_0$ such that for each $k \geq k_0$, the set $\mathcal{G}^k(M)$ of simple $C^k(M)$ metrics in $M$ for which $I_g$ is s-injective is open and dense in the $C^k(M)$ topology. Moreover, for any $g \in \mathcal{G}^k$,*

$$\|f_M^s\|_{L^2(M)} \leq C \|N_g f\|_{\widetilde{H}^2(M_1)}, \quad \forall f \in H^1(M), \tag{1.2}$$

*with a constant $C > 0$ that can be chosen locally uniform in $\mathcal{G}^k$ in the $C^k(M)$ topology.*

Of course, $\mathcal{G}^k$ includes all real analytic simple metrics in $M$, according to Theorem 1.

The analysis of $I_g$ can also be carried out for symmetric tensors of any order, see e.g. [Sh1]. Since we are motivated by the boundary rigidity problem, and to simplify the exposition, we study only tensors of order 2.

Theorem 2 and especially estimate (1.2) allow us to prove in section 5 the following local generic uniqueness result for the non-linear boundary rigidity problem.

**Theorem 3** *Let $k_0$ and $\mathcal{G}^k(M)$ be as in Theorem 2. There exists $k \geq k_0$, such that for any $g_0 \in \mathcal{G}^k$, there is $\varepsilon > 0$, such that for any two metrics $g_1$, $g_2$ with $\|g_m - g_0\|_{C^k(M)} \leq \varepsilon$, $m = 1, 2$, we have the following:*

$$\rho_{g_1} = \rho_{g_2} \text{ on } (\partial M)^2 \text{ implies } g_2 = \psi_* g_1 \tag{1.3}$$

*with some $C^{k+1}(M)$-diffeomorphism $\psi : M \to M$ fixing the boundary.*

Theorem 3 implies the following generic global uniqueness result for simple metrics.

**Theorem 4** *There exists an open dense set $\mathcal{U}$ of pairs of simple metrics in $C^k(M) \times C^k(M)$, with $k$ as in Theorem 3, such that any pair in that set has the property (1.3). Moreover, for any simple $g_1 \in C^k(M)$, the set of simple metrics with the same boundary distance function as $g_1$ on $(\partial M)^2$ is a complement of a dense open set in the $C^k(M)$ topology.*

Finally, we prove a conditional stability estimate of Hölder type. A similar estimate near the Euclidean metric was proven in [W] based on the approach in [SU2].

**Theorem 5** *Let $k_0$ and $\mathcal{G}^k(M)$ be as in Theorem 2. Then for any $\mu < 1$, there exits $k \geq k_0$ such that for any $g_0 \in \mathcal{G}^k$, there is an $\varepsilon_0 > 0$ and $C > 0$ with the property that that for any two metrics $g_1$, $g_2$ with $\|g_m - g_0\|_{C(M)} \leq \varepsilon_0$, and $\|g_m\|_{C^k(M)} \leq A$, $m = 1, 2$, with some $A > 0$, we have the following stability estimate*

$$\|g_2 - \psi_* g_1\|_{C^2(M)} \leq C(A) \|\rho_{g_1} - \rho_{g_2}\|^\mu_{C(\partial M \times \partial M)}$$

*with some diffeomorphism $\psi : M \to M$ fixing the boundary.*

One can formulate a global generic stability result in the spirit of Theorem 4, with a constant uniform in any compact subset of $\mathcal{U}$.

A brief outline of the paper is as follows. In section 3 we construct the analytic parametrix for the normal operator $N_g$. In section 4 we prove Theorem 1. In section 5 we use the results of the previous section to prove Theorem 2, i.e., generic s-injectivity of $I_g$. Theorem 3 about generic uniqueness for the non-linear problem is proved in section 5. The stability estimate of Theorem 5 is proved in section 7. Even though Theorem 5 implies Theorem 3, we still include a proof of the latter in section 6 for convenience of the reader, since it is much shorter.

We would like to thank C. Croke for very useful comments on a previous version of the paper.

## 2 Preliminaries

We start with some basic facts about symmetric 2-tensors on Riemannian manifolds. For more details, we refer to [Sh1].

As we mentioned in the introduction it is enough to prove the results for the case that the manifold is the closure of a bounded open set with smooth boundary which we will denote by $\Omega$.

Let $g \in C^2(\bar\Omega)$ be a Riemannian metric in $\Omega$. We work with symmetric 2-tensors $f = \{f_{ij}\}$ and with 1-tenors (differential forms) $v_j$. We use freely the Einstein summation convention and the convention for raising and lowering indices. We think of $f_{ij}$ and $f^{ij} = f_{kl}g^{ki}g^{lj}$ as different representations of the same tensor. If $(x, \xi)$ is a covector, then its components are denoted by $x^j$, and $\xi_j$, while $\xi^j$ is defined as $\xi^j = \xi_i g^{ij}$. Next, we denote $|\xi|^2 = \xi_i \xi^i$.

The geodesics of $g$ can be also viewed as the $x$-projections of the bicharacteristics of the Hamiltonian $H_g(x, \xi) = \frac{1}{2} g^{ij}(x) \xi_i \xi_j$. The energy level $H_g = 1/2$ corresponds to parameterization with arc-length parameter. For any geodesic $\gamma$, we have $f_{ij}(\gamma(x))\dot\gamma^i(t)\dot\gamma^j(t) = f^{ij}(x)\xi_i \xi_j$, where $(x, \xi) = (x(t), \xi(t))$ is the bicharacteristic with $x$-projection equal to $\gamma$.

Set

$$\Gamma_- := \left\{ (x, \omega) \in T^*\Omega;\ x \in \partial\Omega,\ |\omega| = 1,\ \omega_i v^i < 0 \right\},$$

where $v(x)$ is the outer unit normal to $\partial\Omega$. Introduce the measure

$$d\mu(x, \omega) = |\omega \cdot v(x)|\, dS_x\, dS_\omega \quad \text{on } \Gamma_-,$$

where $dS_x$ and $dS_\omega$ are the surface measures on $\partial\Omega$ and $\{\omega \in T^*_x\Omega;\ |\omega| = 1\}$ in the metric, respectively. If $\partial\Omega$ is given locally by $x^n = 0$, then $dS_x = (\det g)^{1/2} dx^1 \ldots dx^{n-1}$, and $dS_\omega = (\det g)^{-1/2} dS_{\omega_0}$, where $dS_{\omega_0}$ is the Euclidean measure on $S^{n-1}$. Let $(x(t; z, \omega), \xi(t; z, \omega))$ be the bicharacteristic issued from $(z, \omega) \in \Gamma_-$ extended for $t \geq 0$ until the first component reaches $\partial\Omega$ again. Then we set

$$If(z, \omega) = \int f^{ij}(x(t; z, \omega)) \xi_i(t; z, \omega) \xi_j(t; z, \omega)\, dt, \quad (z, \omega) \in \Gamma_-.$$

We will drop the subscript $g$ in $I_g$ unless we have more than one metric and the subscript is needed. By identifying $T\Omega$ with $T^*\Omega$, as explained above, we can think of $\Gamma_-$ as a subset of $T\Omega$, and then

$$If(x,\theta) = \int f_{ij}(\gamma_{x,\theta}(t))\dot{\gamma}^i_{x,\theta}(t)\dot{\gamma}^j_{x,\theta}(t)\,\mathrm{d}t, \quad (x,\theta) \in \Gamma_-,$$

where $\gamma_{x,\theta}$ is the geodesic issued from $x$ in the direction $\theta$. Since $|\theta| = 1$, $t$ is the arc-length parameter. Clearly, $I : C^\infty(\bar{\Omega}) \to L^2(\Gamma_-, \mathrm{d}\mu)$. We define the $L^2$ space of symmetric tensors $f = \{f_{ij}\}$ with inner product

$$(f, h) = \int_\Omega f_{ij} \bar{h}^{ij} (\det g)^{1/2}\,\mathrm{d}x.$$

Similarly, we define the $L^2$ space of 1-tensors (vector fields, that we identify with 1-forms) and the $L^2$ space of functions in $\Omega$. Also, we will work in Sobolev $H^s$ spaces of 2-tensors, 1-forms and functions. In order to keep the notation simple, we will use the same notation $L^2$ (or $H^s$) for all those spaces and it will be clear from the context which one we mean.

The map $I : L^2(\Omega) \to L^2(\Gamma_-, \mathrm{d}\mu)$ is bounded [SU3], and therefore the *normal* operator $N := I^*I$ is a well defined bounded operator in $L^2(\Omega)$. In [SU3], we found that

$$(Nf)_{kl}(x) = \frac{2}{\sqrt{\det g}} \int \frac{f^{ij}(y)}{\rho(x,y)^{n-1}} \frac{\partial \rho}{\partial y^i} \frac{\partial \rho}{\partial y^j} \frac{\partial \rho}{\partial x^k} \frac{\partial \rho}{\partial x^l} \det \frac{\partial^2 (\rho^2/2)}{\partial x \partial y}\,\mathrm{d}y, \quad x \in \Omega. \tag{2.1}$$

We always assume that $g$ is extended as a simple metric in a small neighborhood of $\Omega$. Also, we always extend functions or tensors defined in $\Omega$, or similar domains, as 0 outside the domain. Then $Nf$ is well defined near $\Omega$ as well and is independent of small perturbations $\partial\Omega_1$ of $\partial\Omega$, as long as $\Omega_1 \supset \Omega$.

Given a symmetric 2-tensor $f = f_{ij}$, we define the 1-tensor $\delta f$ called *divergence* of $f$ by

$$[\delta f]_i = g^{jk} \nabla_k f_{ij},$$

where $\nabla_k$ are the covariant derivatives of the tensor $f$. Given a 1-tensor (vector field) $v$, we denote by $dv$ the 2-tensor called symmetric differential of $v$:

$$[dv]_{ij} = \frac{1}{2}\left(\nabla_i v_j + \nabla_j v_i\right).$$

Operators $d^s$ and $-\delta^s$ are formally adjoint to each other in $L^2(\Omega)$. It is easy to see that for each smooth $v$ with $v = 0$ on $\partial\Omega$, we have $I(dv) = 0$. This follows from the identity

$$\frac{\mathrm{d}}{\mathrm{d}t} v_i(\gamma(t))\dot{\gamma}^i(t) = [dv(\gamma(t))]_{ij}\dot{\gamma}^i(t)\dot{\gamma}^j(t). \tag{2.2}$$

It is known that for $g$ smooth enough (see Lemma 1 below), each symmetric tensor $f \in L^2(\Omega)$ admits unique orthogonal decomposition $f = f^s + dv$ into a *solenoidal* tensor $\mathcal{S}f := f^s$ and a *potential* tensor $\mathcal{P}f := dv$, such that both terms are in $L^2(\Omega)$, $f^s$ is solenoidal, i.e., $\delta f^s = 0$ in $\Omega$, and $v \in H^1_0(\Omega)$ (i.e., $v = 0$ on $\partial\Omega$). In order to construct this decomposition, introduce the operator $\Delta^s = \delta d$ acting on vector fields. This operator is elliptic in $\Omega$, and the Dirichlet problem satisfies the Lopatinskii condition. Denote by $\Delta^s_D$ the Dirichlet realization of $\Delta^s$ in $\Omega$. Then

$$v = \left(\Delta^s_D\right)^{-1} \delta f, \quad f^s = f - d\left(\Delta^s_D\right)^{-1} \delta f. \tag{2.3}$$

Notice that even when $f$ is smooth and $f = 0$ on $\partial\Omega$, then $f^s$ does not need to vanish on $\partial\Omega$. To stress on the dependence on the domain, when needed, we will use the notation $v_\Omega$ and $f^s_\Omega$ as well.

Operators $\mathcal{S}$ and $\mathcal{P}$ are orthogonal projectors, see also Lemma 1 below. The problem about the s-injectivity of $I$ then can be posed as follows: if $If = 0$, show that $f^s = 0$, in other words, show that $I$ is injective on the subspace $\mathcal{S}L^2$ of solenoidal tensors. Note that $N\mathcal{S} = \mathcal{S}N, \mathcal{P}N = N\mathcal{P} = 0$.

In [SU3], we analyzed $N$. We showed that for simple metrics, it is a $\Psi$DO of order $-1$, as a consequence of (2.1). The operator $N$ is not elliptic, it actually vanishes on $\mathcal{P}L^2(\Omega)$, and its principal symbol $\sigma_p(N)_{ijkl}$ vanishes on the Fourier transforms of locally potential (symmetric) tensors. On the complement of this linear space however, $\sigma_p(N)$ is elliptic. This allows us to construct a parametrix, and we will use this construction in next two sections. To obtain $f^s$ modulo smoothing operators, however, requires additional efforts, see Proposition 4, because we work in a domain with boundary and $f^s$ is defined by solving a BVP as well, see (2.3). Nevertheless, those arguments allowed us to prove the following a priori stability estimate for the linear problem [SU3] for simple smooth metrics $g$:

$$\|f^s_\Omega\|_{L^2(\Omega)} \leq C\|Nf\|_{\widetilde{H}^2(\Omega_1)} + C_s\|f\|_{H^{-s}(\Omega_1)}, \quad \forall f \in H^1(\Omega), \quad \forall s > 0. \tag{2.4}$$

The norm $\|\cdot\|_{\widetilde{H}^2}$ is introduced in (5.2) and $\Omega_1 \supset \Omega$ is a small neighborhood of $\Omega$, where $g$ is still simple. If $I$ is s-injective, then we showed that one can choose $C_s = 0$. Notice that there is a loss of one derivative in (2.4).

In our exposition, we will often use boundary normal (semi-geodesic) coordinates $(x', x^n)$ near a boundary point. They are defined such that $x^n = 0$ defines $\partial\Omega$, $x^n > 0$ in $\Omega$, and $x^n = \rho(x, \partial\Omega)$. The metric $g$ in those coordinates satisfies $g_{in} = \delta_{in}$, and $\Gamma^i_{nn} = \Gamma^n_{in} = 0$, $\forall i$. We also use the convention that all Greek indices take values from 1 to $n-1$.

At the end of this section we prove a lemma that justifies the decomposition of a symmetric $L^2$–tensor into a solenoidal and a potential part.

**Lemma 1** *For each $g \in C^1(\bar{\Omega})$, the maps*

$$(\Delta^s_D)^{-1} : H^{-1}(\Omega) \to H^1_0(\Omega), \quad \mathcal{P}, \mathcal{S} : L^2(\Omega) \longrightarrow L^2(\Omega)$$

*are bounded and depend continuously on $g$.*

*Proof:* We consider $(\Delta^s_D)^{-1}$ first (see also [N, Lemma 4.2], where $\Delta^s$ is replaced by $\nabla \cdot \gamma(x)\nabla$ and $\gamma$ is assumed to be uniformly Lipschitz).

We start with the fact that for a fixed smooth $g_0$, $(\Delta^s_{g_0,D})^{-1}$ is bounded in the spaces above [Ta, p. 307]. Let $g_0$ and $g$ be two smooth metrics, and let $u$ and $v$ be in $H^{-1}(\Omega)$. Then

$$(\Delta^s_{g,D})^{-1} - (\Delta^s_{g_0,D})^{-1} = (\Delta^s_{g,D})^{-1} \left(\Delta^s_{g_0,D} - \Delta^s_{g,D}\right) (\Delta^s_{g_0,D})^{-1} \tag{2.5}$$

To estimate the difference in the r.h.s. above with the regularity assumptions on $g$, choose $u$ and $v$ in $H^1_0(\Omega)$. Then

$$\left|\langle(\Delta^s_{g,D} - \Delta^s_{g_0,D})u, v\rangle\right| = \left|\langle d_{g_0}u, d_{g_0}v\rangle - \langle d_g u, d_g v\rangle\right|$$
$$\leq C\|g - g_0\|_{C^1}(\|g\|_{C^1} + \|g_0\|_{C^1})\|u\|_{H^1}\|v\|_{H^1}.$$

Let $\|g - g_0\|_{C^1} \leq \varepsilon$. Then for $\varepsilon \ll 1$, $\|g - g_0\|_{C^1}(\|g\|_{C^1} + \|g_0\|_{C^1}) \leq C(g_0)\varepsilon$. As a consequence,

$$\|(\Delta^s_{g,D})^{-1}\|_{H^{-1} \to H^1} \leq \|(\Delta^s_{g_0,D})^{-1}\|_{H^{-1} \to H^1} \left(1 + C(g_0)\varepsilon\|(\Delta^s_{g,D})^{-1}\|_{H^{-1} \to H^1}\right).$$

This implies that for $\varepsilon \ll 1$, the l.h.s. above is uniformly bounded by a constant depending on $g_0$.

$$\|(\Delta_{g,D}^s)^{-1}\|_{H^{-1}\to H^1} \leq (1-C\varepsilon)^{-1}\|(\Delta_{g_0,D}^s)^{-1}\|_{H^{-1}\to H^1},$$

provided that $\|g-g_0\|_{C^3} \leq \varepsilon < 1/C$ with $C$ depending on $\Omega$ only. Going back to (2.5), we conclude

$$\|(\Delta_{g,D}^s)^{-1} - (\Delta_{g_0,D}^s)^{-1}\|_{H^{-1}\to H^1} \leq C\|g-g_0\|_{C^1} \tag{2.6}$$

with $C$ a uniform constant in any small enough $C^1$-neighborhood of a fixed $g_0$. This inequality allows us to define $(\Delta_{g,D}^s)^{-1} : H^{-1}(\Omega) \to H_0^1(\Omega)$ for any metric $g \in C^1(\bar{\Omega})$ by approximating with smooth $g$. Moreover, we get that the resolvent above is continuous in $g$ and (2.6) still holds. As a consequence, $\mathcal{S}$ and $\mathcal{P}$ are also continuous in $g \in C^1$ as operators in $L^2(\Omega)$. This completes the proof of the lemma. $\square$

Remark also that the lemma above admits the following easy generalization: for $s = 0, 1, \ldots$, the resolvent in the lemma also continuously maps $H^{s-1}$ into $H^{s+1} \cap H_0^1$, similarly, $\mathcal{P}$ and $\mathcal{S}$ are bounded in $H^s$, if $g \in C^k$, $k \gg 1$ (depending on $s$). Moreover those operators depend continuously on $g$.

## 3 The analytic parametrix

In what follows, "analytic" always means real analytic.

Assume that $g$ is a simple analytic metric in $\bar{\Omega}$. Our goal is to reconstruct $f_{\Omega_1}^s$ from $Nf$ up to an analytic-regularizing operator, where $\Omega_1$ is a slightly larger domain. This is the key step towards proving Theorem 1 in next section.

We are going to use the analytic $\Psi$DO calculus, see [Tre]. Analytic $\Psi$DO have been used in integral geometry before, see e.g., [BQ] for uniqueness results for the Euclidean Radon transform with analytic weights.

Next, we will follow the parametrix construction in [SU3] in the new situation, where $g$ is analytic. Since $N$ is not elliptic, we modify it to get an elliptic operator of order zero first:

$$W = \chi N + N_0 \mathcal{P}_{\Omega_2}, \tag{3.1}$$

where $\chi$ is a smooth cutoff function equal to 1 in a neighborhood of $\bar{\Omega}_1$, supported in a larger neighborhood of $\bar{\Omega}_1$ where $g$ is still simple; and $N_0$ is the operator with symbol $(\xi_1^2 + \ldots + \xi_n^2)^{-1/2}$, i.e., an integral operator with kernel $c_n|x-y|^{-n+1}$. Recall that $\mathcal{P}_{\Omega_2} = d(\Delta_{\Omega_2,D}^s)^{-1}\delta$. Here $\Omega \subset \Omega_1 \subset \Omega_2$, and $\Omega_1$ is a small strictly convex neighborhood of $\Omega$ with analytic boundary, $\Omega_2$ is related to $\Omega_1$ in the same way, and we extend $g$ analytically near $\Omega_2$. Inside $\Omega_2$, and therefore, on $\bar{\Omega}_1$, the operator $W$ is an elliptic $\Psi$DO [SU3].

Similarly to [SU3], we have.

**Lemma 2** *There exists $\delta > 0$ such that in $U = \{(x,y) \in \Omega_2 \times \Omega_2; |x-y| < \delta\}$ we have*

$$\rho^2(x,y) = G_{ij}^{(1)}(x,y)(x-y)^i(x-y)^j,$$
$$\frac{\partial \rho^2(x,y)}{\partial x^j} = 2G_{ij}^{(2)}(x,y)(x-y)^i,$$
$$\frac{\partial^2 \rho^2(x,y)}{\partial x^i \partial y^j} = 2G_{ij}^{(3)}(x,y),$$

where $G_{ij}^{(1)}$, $G_{ij}^{(2)}$, $G_{ij}^{(3)}$ are analytic in $U$, positively defined, and we have

$$G_{ij}^{(1)}(x,x) = G_{ij}^{(2)}(x,x) = G_{ij}^{(3)}(x,x) = g_{ij}(x).$$

*Proof of Lemma 2:* Let the covector $\xi$ be defined as $\xi = \xi(x,y) = \exp_x^{-1} y$. Then $\xi(x,x) = 0$, therefore

$$\xi_i = A_{ij}(x,y)(x^j - y^j) \quad \text{with} \quad A_{ij}(x,y) = \int_0^1 \partial_{y_j} \xi_i(x, x + t(y-x)) \, dt. \tag{3.2}$$

The latter is a well defined analytic function for $x - y$ small enough since then the line segment $[x, y]$ along which we integrate does not leave $\Omega_2$. For $x$ and $y$ far apart, it may leave $\Omega_2$ which is geodesically convex but not necessarily convex w.r.t. the Euclidean metric. It is easily seen that $A_{ij}(x,x) = g_{ij}(x)$, and $g^{ij}(x)\xi_i\xi_j = \rho^2(x,y)$, so (3.2) implies the lemma. $\square$

**Proposition 1** *N and W are analytic $\Psi$DOs in $\Omega_2$.*

*Proof:* We analyze $N$ first. Recall (2.1). Let $V$ be open such that $V \times V \subset U$, supp $f \subset V$. Then for $x \in V$,

$$[Nf]_{ij}(x) = \int \check{M}_{ijkl}(x, y, x - y) f^{kl}(y) \, dy$$

with

$$M_{ijkl}(x, y, \xi) = \int e^{-i\xi \cdot z} \left(G^{(1)} z \cdot z\right)^{\frac{-n+1}{2} - 2} \tag{3.3}$$

$$\times [G^{(2)} z]_i [G^{(2)} z]_j [\widetilde{G}^{(2)} z]_k [\widetilde{G}^{(2)} z]_l \frac{\det G^{(3)}}{\sqrt{\det g}} \, dz,$$

and $\widetilde{G}_{ij}^{(2)}(x,y) = G_{ij}^{(2)}(y,x)$. It is convenient to make the change $z' = \left(G^{(1)}(x,y)\right)^{1/2} z$ above to get $M(x, y, \xi) = \widetilde{M}(x, y, (G^{(1)}(x,y))^{-1/2} \xi)$, where

$$\widetilde{M}_{ijkl}(x, y, \xi) = \int e^{-i\xi \cdot z} |z|^{-n-3} [G^{(2)}(G^{(1)})^{-1/2} z]_i [G^{(2)}(G^{(1)})^{-1/2} z]_j \tag{3.4}$$

$$\times [\widetilde{G}^{(2)}(G^{(1)})^{-1/2} z]_k [\widetilde{G}^{(2)}(G^{(1)})^{-1/2} z]_l \det(G^{(1)})^{-n/2} \frac{\det G^{(3)}}{\sqrt{\det g}} \, dz,$$

As a Fourier transform of a (positively) homogeneous in $z$ distribution, $\widetilde{M}$ is homogeneous in $\xi$ of order $-1$. It is analytic function of all variables for $\xi \neq 0$. To prove this, write

$$\widetilde{M}(x, y, \xi) = \int e^{-i\xi \cdot z} |z|^{-n+1} m(x, y, \theta) \, dz, \quad \theta = z/|z|$$

and pass to polar coordinates $z = r\theta$. Since $m$ is an even function of $\theta$, we get (see also [H, Theorem 7.1.24])

$$\widetilde{M}(x, y, \xi) = \pi \int_{|\theta|=1} m(x, y, \theta) \delta(\theta \cdot \xi) \, d\theta,$$

and our claim follows since $m$ is analytic function of all its variables in the integral above.

Let $\chi \in C_0^\infty$. We will prove first that $\chi(\xi) M_{ijkl}(x, y, \xi)$ is an amplitude of an analytic-regularizing operator for $(x, y) \in U$. Indeed,

$$\begin{aligned}(\chi M)(x, y, D) f &= (2\pi)^{-n} \int \int e^{i(x-y)\cdot\xi} \chi(\xi) M(x, y, \xi) f^{kl}(y) \, dy \, d\xi \\ &= (2\pi)^{-n} \int_{S^{n-1}} \int_0^\infty \int e^{i(x-y)\cdot r\theta} \chi(r\theta) M(x, y, \theta) f^{kl}(y) r^{n-2} \, dy \, dr \, d\theta,\end{aligned}$$

and the analyticity follows from this representation.

Next, $(1 - \chi(\xi)) M_{ijkl}(x, y, \xi)$ is an analytic amplitude [Tre, Definition V.2.1] for $(x, y) \in U$. The estimates needed to justify this statement follow from the homogeneity of $M$ and the Cauchy integral formula.

The arguments above prove that for any $x_0 \in \Omega_2$, there exists a neighborhood $V_{x_0}$ of $x_0$, such that $M$ is an analytic amplitude for $(x, y) \in V_{x_0} \times V_{x_0}$, therefore $N$ is an analytic $\Psi$DO in $V_{x_0}$. To prove that $N$ is an analytic $\Psi$DO in the whole $\Omega_2$, we will follow the proof of [Tre, Theorem V.3.4]. The statement follows from the fact that the kernel of $N$ is analytic away from the diagonal, which, combined with what we proved above implies easily that $N$ is analytic pseudo-local in the whole $\Omega_2$. More precisely, one can define the analytic formal symbol

$$\exp\{\partial_\xi D_y\} M(x, y, \xi)|_{y=x}$$

and this symbol defines an equivalence class of analytic $\Psi$DOs in a neighborhood of $\Omega_2$. One can build a true pseudo-analytic symbol $\tilde{a}(x, \xi)$ in $\Omega_2$ based on the formal series above as in [Tre]. For any sufficiently small open set $V$, one has that $(N - a(x, D))u$ is analytic in $\Omega_2$ for $u \in \mathcal{E}'(V)$, and one can easily extend this to any distribution supported in $\Omega_2$. This completes the proof for $N$.

Consider next $(\Delta^s_{\Omega_2, D})^{-1}$. The operator $\Delta^s$ is an analytic elliptic $\Psi$DO, therefore, it has a parametrix $P$, that is analytic $\Psi$DO in $\Omega_2$, such that $\Delta^s P$ is analytic-regularizing in $\bar{\Omega}_2$ (we need to work in a bit larger domain in order to cover $\bar{\Omega}_2$). Let $u = (\Delta^s_{\Omega_2, D})^{-1} f$. If $\operatorname{supp} f \subset \Omega_2$, then $u - Pf$ solves equation of the kind (3.5) below with analytic coefficients, therefore, by the interior analytic regularity, $u - Pf$ is analytic in $\Omega_2$. This shows that $(\Delta^s_{\Omega_2, D})^{-1}$ equals $P$ up to an analytic-regularizing operator in any compact subset, therefore, $(\Delta^s_{\Omega_2, D})^{-1}$ is analytic $\Psi$DO in $\Omega_2$.

The remaining operators in (3.1) are clearly analytic $\Psi$DOs. This completes the proof of Proposition 1. □

Next step is to reconstruct $f^s_{\Omega_1}$ from $Nf$ up to analytic function. We need first a lemma about analyticity up to the boundary of solutions of $\Delta^s v = u$:

**Lemma 3** *Let $x_0 \in \partial\Omega$, and assume that the metric $g$, and the vector fields $u$, $v_0$ are analytic in a (two-sided) neighborhood of $x_0$, and that $\partial\Omega$ is analytic near $x_0$. Let the vector field $v$ solve*

$$\Delta^s v = u \quad \text{in } \Omega, \quad v|_{\partial\Omega} = v_0. \tag{3.5}$$

*Then $v$ extends as analytic function in some (two-sided) neighborhood of $x_0$.*

*Proof:* The lemma follows directly from [MN]. One can first extend $v_0$ near $x_0$ as analytic function, and subtract from $v$ certain function analytic near $x_0$ that reduces the problem to one with $v_0 = 0$. Next, we observe that the principal symbol of $-\Delta^s$ is a positive matrix for $\xi \neq 0$, hence the system above is strongly elliptic in the terminology of [MN], and the result follows (see also [Tre]). □

**Proposition 2** *There exists a bounded operator $P : H^1(\Omega_1) \to L^2(\Omega_1)$, such that for any symmetric 2-tensor $f \in L^2(\Omega)$ we have*

$$f^s_{\Omega_1} = PNf + Kf,$$

*with $Kf$ analytic in $\bar{\Omega}_1$. Moreover, $P$ is an analytic $\Psi DO$ in a neighborhood of $\Omega_1$ of order 1.*

*Proof:* We follow the proof of Theorem 2 in [SU3], where $g$ is smooth only. First, we construct a parametrix $L$ of $W$ in $\Omega_2$, see [Tre, Theorem V.3.3]. There exists $L : \mathcal{D}'(\Omega_2) \to \mathcal{E}'(\Omega_2)$ such that $L$ is an analytic $\Psi DO$ of order 1 in a neighborhood of $\bar{\Omega}_1$, and such that $LW = \mathrm{Id} + K_1$ near $\bar{\Omega}_1$, where $K_1 f \in \mathcal{A}(\bar{\Omega}_1)$ for any $f \in L^2(\Omega_1)$. Chose a smooth $\chi$ supported in $\Omega_2$ such that $\chi = 1$ near $\bar{\Omega}_1$. Then the equality above implies $\chi L W \chi = \chi^2 + \chi K_1 \chi$. Apply $\mathcal{S}_{\Omega_2}$ to the left and right to get

$$\mathcal{S}_{\Omega_2} L W \mathcal{S}_{\Omega_2} = \mathcal{S}_{\Omega_2} + K_2$$

with $K_2$ having the property that $K_2 f \in \mathcal{A}(\bar{\Omega}_1)$ for any $f \in L^2(\Omega_1)$. To prove the latter, we use the analytic pseudolocal property of analytic $\Psi DOs$.

We have $W \mathcal{S}_{\Omega_2} = N$. Therefore, setting $P = \mathcal{S}_{\Omega_2} L$, we get

$$PN = \mathcal{S}_{\Omega_2} + K_2 \quad \text{in } \Omega_2. \tag{3.6}$$

Note that we have showed that $K_2$ maps $L^2(\Omega_1)$ into $\mathcal{A}(\bar{\Omega}_1)$, but not into $\mathcal{A}(\bar{\Omega}_2)$.

Next, compare $f^s_{\Omega_1}$ and $f^s_{\Omega_2}$ for $f \in L^2(\Omega)$. We have $f^s_{\Omega_i} = f - dv_{\Omega_i}$, where $v_{\Omega_i} = (\Delta^s_{\Omega_i,D})^{-1} \delta f$, $i = 1, 2$. Thus $f^s_{\Omega_1} = f^s_{\Omega_2} + dw$ in $\Omega_1$, where the vector field $w = v_{\Omega_2} - v_{\Omega_1} \in H^1(\Omega_1)$ solves

$$\Delta^s w = 0 \quad \text{in } \Omega_1, \quad w|_{\partial \Omega_1} = v_{\Omega_2}. \tag{3.7}$$

Since supp $f$ is disjoint from $\partial \Omega_1$, we get $v_{\Omega_2} \in \mathcal{A}(\partial \Omega_1)$. By Lemma 3, $w \in \mathcal{A}(\bar{\Omega}_1)$, thus $f \mapsto dw|_{\Omega_1}$ is a linear operator mapping $L^2(\Omega)$ into $\mathcal{A}(\bar{\Omega}_1)$. Then the relation

$$f^s_{\Omega_1} = f^s_{\Omega_2} + dw = PNf - K_2 f + dw$$

completes the proof of the proposition. □

## 4  S-injectivity for analytic metrics; proof of Theorem 1

In this section, we prove Theorem 1. We start with a recovery at the boundary result. Next lemma generalizes Lemma 2.3 in [Sh2] by proving that actually $v = 0$ with all derivatives at $\partial \Omega$. On the other hand, in [Sh2], $\Omega$ does not need to be convex. Also, the lemma below can be considered as a linear version of Theorem 2.1 in [LSU]: if two metrics have the same boundary distance function, then in boundary normal coordinates, they have the same derivatives of all orders at $\partial \Omega$.

**Lemma 4** *Let $g$ be a smooth, simple metric in $\Omega$ and let $f_{ij}$ be a smooth symmetric tensor. Assume that $If = 0$. Then there exists a smooth vector field $v$ vanishing on $\partial \Omega$, supported near $\partial \Omega$, such that for $\tilde{f} = f - dv$ we have $\partial^m \tilde{f}|_{\partial \Omega} = 0$ for any multiindex $m$.*

*Moreover, if $g$ and $f$ are analytic in a (two-sided) neighborhood of $\partial \Omega$, then $v$ can be chosen so that $\tilde{f} = 0$ near $\partial \Omega$ as well.*

*Proof:* We fix $x_0 \in \partial\Omega$ and below we work in some neighborhood of $x_0$. Assume that $x = (x', x^n)$ are boundary normal coordinates near $x_0$, therefore there we have $g_{in} = \delta_{in}$, $\forall i$. We will find a vector filed $v$ vanishing on $\partial\Omega$ such that for $\tilde{f} := f - dv$ we have $\tilde{f}_{in} = 0$ for $i = 1, 2, \ldots, n$. The latter is equivalent to

$$\nabla_n v_i + \nabla_i v_n = 2 f_{in}, \quad v|_{x^n=0} = 0, \quad i = 1, \ldots, n. \tag{4.1}$$

Recall that $\nabla_i v_j = \partial_i v_j - \Gamma_{ij}^k v_k$, and that in those coordinates, $\Gamma_{nn}^k = \Gamma_{kn}^n = 0$. We solve (4.1) for $i = n$ first by integration, then $\nabla_n v_n = \partial_n v_n = f_{nn}$; this gives us $v_n$. Next, we solve the remaining linear system of $n - 1$ equations for $i = 1, \ldots, n-1$ that is of the form $\nabla_n v_i = 2 f_{in} - \nabla_i v_n$, or, equivalently,

$$\partial_n v_i - 2\Gamma_{ni}^\alpha v_\alpha = 2 f_{in} - \partial_i v_n, \quad v_i|_{x^n=0} = 0, \quad i = 1, \ldots, n-1. \tag{4.2}$$

(recall that $\alpha = 1, \ldots, n-1$). Clearly, if $g$ and $f$ are analytic near $\partial\Omega$, then so is $v$.

We have $I\tilde{f} = 0$ for $(x, \xi)$ such that $x \in \partial\Omega$ is close to $x_0$, $|\xi| = 1$, and its normal component is small enough. This guarantees that the geodesic $\gamma_{x,\xi}$ issued from $(x, \xi)$ hits the boundary again at a point close to $x_0$, where $v = 0$. We can adapt the proof of Theorem 2.1 in [LSU] to our situation. For the sake of completeness, we will repeat those arguments. It is enough to prove that

$$\partial_n^j \tilde{f}_{\alpha\beta}|_{x=x_0} = 0, \quad \forall j = 0, 1, \ldots, \quad \forall \alpha, \beta = 1, \ldots, n-1. \tag{4.3}$$

Indeed, if (4.3) holds, then in the same way we prove (4.3) for $x \in \partial\Omega$ close to $x_0$, so we can differentiate (4.3) w.r.t. $x'$ to we get that all derivatives of $\tilde{f}$ on $\partial\Omega$ vanish.

Notice that (4.3) is obvious for $j = 0$. Assume that there is $j \geq 1$ such that (4.3) is not true. The Taylor expansion of $\tilde{f}$ then implies that $\exists \xi_0$ of unit length tangent to $\partial\Omega$ such that $\tilde{f}_{\alpha\beta}(x)\xi^\alpha\xi^\beta$ is either (strictly) positive or negative for $x^n > 0$ and $x'$ both sufficiently small and $\xi$ close to $\xi_0$. Notice that in the summation above, we have $\alpha < n$ and $\beta < n$ because $\tilde{f}_{in} = \tilde{f}_{ni} = 0$. Therefore, $I\tilde{f}$ is either (strictly) positive or negative for all $(x, \xi) \in \Gamma_-$ close enough to $(x_0, \xi_0)$ and this is a contradiction.

To make that construction global near $\partial\Omega$, it is enough to note that equation (4.1) is invariant under coordinate changes, so the local construction in fact yields a global one near $\partial\Omega$, see also [Sh2, Lemma 2.2]. Finally, we cut $v$ near $\partial\Omega$ to complete the proof.

If $g$ is analytic up to $\partial\Omega$, then as pointed out above, $v$ is analytic near $\partial\Omega$, up to $\partial\Omega$. Therefore, the same is true for $\tilde{f}$, thus $\tilde{f} = 0$ near $\partial\Omega$. $\square$

Next, we introduce global semi-geodesic coordinates in $\bar{\Omega}$ already used in [SU2], [SU3].

**Lemma 5** *Let $g \in C^k(\bar{\Omega})$, $k \geq 2$, be a simple metric in $\Omega$. Then there exists a $C^{k-1}$ diffeomorphism $\psi : \bar{\Omega} \to \psi(\bar{\Omega})$, such that in the coordinates $y = \psi(x)$, the metric $g$ has the property*

$$g_{in} = \delta_{in}, \quad i = 1, \ldots, n. \tag{4.4}$$

*Moreover, if $g \in \mathcal{A}(\bar{\Omega})$, then $\psi \in \mathcal{A}(\bar{\Omega})$.*

*Proof:* The proof is essentially given in [SU3] and is based on defining the so-called normal coordinates near a fixed point. Let $\Omega_1 \supset\supset \Omega$ be as above and fix $x_0 \in \partial\Omega_1$. Then $\exp_{x_0}^{-1} : \Omega_1 \to \exp_{x_0}^{-1}(\Omega_1)$ is a diffeomorphism by our simplicity hypothesis. Choose a Cartesian coordinate system $\xi$ in the tangent space, so that $\xi^n = 0$ is tangent to the boundary of $\exp_{x_0}^{-1}(\Omega_1)$ at $\xi = 0$. Introduce polar coordinates $\xi = r\theta$ in $\exp_{x_0}^{-1}(\Omega_1)$, where $g_{ij}(x_0)\theta^i\theta^j = 1$, $r > 0$. By the strong convexity assumption, $\theta^n$ has a positive lower bound in a neighborhood of the closure of $\exp_{x_0}^{-1}(\Omega_1)$, the same is true for $r$. Then we set $y' = \theta'/\theta_n$, $y_n = r$.

The spheres $r = $ const. are orthogonal to the geodesics $\theta = $ const. by the Gauss lemma. Moreover, $r$ is the arc-length along those geodesics. Passing to the $y$-coordinates, we get that the lines $y' = $ const. are geodesics orthogonal to the planes $y^n = $ const., with $y^n$ arc-length parameter, This proves (4.4).

Clearly, if $g \in \mathcal{A}(\bar{\Omega})$, then the coordinate change above is analytic as well. □

Lemma 5 allows us to assume, without loss of generality, that $g$ satisfies (4.4).

*Proof of Theorem 1:* We work in the semi-geodesic coordinates above. Assume that $g \in \mathcal{A}(\bar{\Omega})$, and let $f \in L^2(\Omega)$ be such that $If = 0$. Then, by Proposition 2, $f_{\Omega_1}^s \in \mathcal{A}(\bar{\Omega}_1)$. Clearly, $If_{\Omega_1}^s = 0$ as well.

Let $v_1$ be the $v$ in Lemma 4, so that $\tilde{f} := f_{\Omega_1}^s - dv_1$ vanishes near $\partial\Omega_1$. Similarly to (4.1) (but now the coordinates are different), we solve

$$\nabla_n v_i + \nabla_i v_n = 2\tilde{f}_{in}, \quad v|_{(\partial\Omega_1)_-} = 0, \qquad (4.5)$$

where $\partial\Omega_\pm$ is the set of all boundary points $x$ for which $(x, e_n) \in \Gamma_\pm$. As before, we first determine $v_n$ by integrating $\partial_n v_n = \tilde{f}_{nn}$ and taking into account the zero boundary condition. Then $v_n = 0$ in a neighborhood $U$ of $\overline{(\partial\Omega_1)_-}$. Next, we solve the remaining linear system (4.2) along the lines parallel to $e_n$ with boundary conditions as in (4.5). We get that $v_i$, $i = 1, \ldots, n-1$ vanish in $U$ as well. For $f^\sharp = f - dv_1 - dv$ we then have that $f^\sharp = 0$ in $U$, and satisfies $f_{in}^\sharp = 0$, $i = 1, \ldots, n$. Moreover, $v_1 + v = 0$ on $(\partial\Omega_1)_-$. On the other hand, there is unique $v^\sharp \in C(\bar{\Omega}_1)$ with the property that for $f^\sharp := f_{\Omega_1}^s - dv^\sharp$ we have $f_{in}^\sharp = 0$, $v^\sharp = 0$ on $(\partial\Omega_1)_-$, and this $v^\sharp$ solves (4.5) with $\tilde{f}$ replaced by $f_{\Omega_1}^s$, so $v^\sharp = v_1 + v$. Since all coefficients in the latter system are analytic, and so is $\partial\Omega_1$, we get that $v^\sharp$ is analytic in $\bar{\Omega} \setminus \partial(\partial\Omega_1)$, i.e., everywhere in $\bar{\Omega}_1$ with a possible exclusion of the points on $\partial\Omega_1$, where $e_n$ is tangent to $\partial\Omega_1$. The same conclusion therefore holds for $f^\sharp$. On the other hand, $f^\sharp = 0$ in $U$, and $U$ includes a neighborhood of $\partial(\partial\Omega_1)$. By analytic continuation, $f^\sharp = 0$ in $\bar{\Omega}_1$.

Thus $f_{\Omega_1}^s = dv^\sharp$ in $\bar{\Omega}_1$, and $v^\sharp = 0$ on $(\partial\Omega_1)_-$. Since we know that $If_{\Omega_1}^s = 0$, by integrating $f_{\Omega_1}^s = dv^\sharp$ along geodesics connecting $(\partial\Omega_1)_-$ and $(\partial\Omega_1)_+$, and using (2.2), we get that $v^\sharp = 0$ in $(\partial\Omega_1)_+$ as well, and by continuity, $v^\sharp = 0$ on the whole $\partial\Omega_1$. This yields $f_{\Omega_1}^s = 0$. Since $\mathrm{supp}\, f \subset \bar{\Omega}$, this easily implies (see next Proposition) that $\mathrm{supp}\, v \subset \bar{\Omega}$, as well.

This concludes the proof of Theorem 1. □

The following elementary statements was used above and it is worth stating separately.

**Proposition 3** *Let $f = dv$, $v|_{\partial\Omega} = 0$, and $v \in C^1(\bar{\Omega})$. Then $\mathrm{supp}\, v$ is included in the smallest simply connected closed set containing $\mathrm{supp}\, f$.*

*Proof:* Let $K$ be the set in the proposition and choose $y \in \Omega \setminus K$. Then there exists a polygon $p = \gamma_1 \cup \gamma_2 \cup \ldots \cup \gamma_m$, each segment $[0,1] \ni t \mapsto \gamma_j$, $j = 1, \ldots, m$ of which is a geodesic, such that $p$ connects some $z \in \partial\Omega$ and $y$. Integrate (2.2) along $\gamma_1$, using the condition $v|_{\partial\Omega} = 0$, to get $v_i \eta^i = 0$ at $\gamma_1(1)$, where $\eta = \dot{\gamma}_1(1)$ is the velocity vector at the endpoint $\gamma_1(1)$ of $\gamma_1$ (different from $z$). By perturbing the initial point $z = \gamma_1(0)$ of $\gamma_1$ a little, and using the simplicity assumption, we get that $v(\gamma_1(1)) = 0$. Similarly, we get that $v = 0$ near $\gamma_1(1)$. Now, we repeat the same argument for $\gamma_2$, etc., until we get $v(y) = 0$. □

# 5   Generic s-injectivity of $I$; proof of Theorem 2

In this section, we prove that the set $\mathcal{G}^k$ is open in the $C^k(\bar{\Omega})$ topology for some $k \gg 1$.

First we recall and modify a little some results in [SU3]. We introduce the norm $\|\cdot\|_{\widetilde{H}^2(\Omega_1)}$ of $Nf$ in $\Omega_1 \supset \Omega$ in the following way. Choose $\chi \in C_0^\infty$ equal to 1 near $\partial\Omega$ and supported in a small neighborhood of $\partial\Omega$ and let $\chi = \sum_{j=1}^J \chi_j$ be a partition of $\chi$ such that for each $j$, on $\operatorname{supp} \chi_j$ we have coordinates $(x'_j, x_j^n)$, with $x_j^n$ a normal coordinate as above. Set

$$\|f\|_{\widetilde{H}^1}^2 = \int \sum_{j=1}^J \chi_j \left( \sum_{i=1}^{n-1} |\partial_{x_j^i} f|^2 + |x_j^n \partial_{x_j^n} f|^2 + |f|^2 \right) dx, \tag{5.1}$$

$$\|Nf\|_{\widetilde{H}^2(\Omega_1)} = \sum_{i=1}^n \|\partial_{x^i} Nf\|_{\widetilde{H}^1} + \|Nf\|_{H^1(\Omega_1)}. \tag{5.2}$$

In other words, in addition to derivatives up to order 1, $\|Nf\|_{\widetilde{H}^2(\Omega_1)}$ includes also second derivatives near $\partial\Omega$ but they are realized as first derivatives of $\nabla Nf$ tangent to $\partial\Omega$.

The reason to use the $\widetilde{H}^2(\Omega_1)$ norm, instead of the stronger $H^2(\Omega_1)$ one, is that this allows us to work with $f \in H^1(\Omega)$, not only with $f \in H_0^1(\Omega)$, since for such $f$, extended as 0 outside $\Omega$, we still have that $Nf \in \widetilde{H}^2(\Omega_1)$, see [SU3].

The following proposition is a modification of the results in section 6 in [SU3].

**Proposition 4** *Let $g \in C^k(\bar{\Omega})$ be simple. Then for any $s = 1, 2, \ldots$, there exists $k > 0$ and a bounded linear operator*

$$Q : \widetilde{H}^2(\Omega_1) \longrightarrow \mathcal{S}L^2(\Omega), \tag{5.3}$$

*such that*

$$QNf = f_\Omega^s + Kf, \quad \forall f \in H^1(\Omega), \tag{5.4}$$

*where $K : H^1(\Omega) \to \mathcal{S}H^{1+s}(\Omega)$ extends to $K : L^2(\Omega) \to \mathcal{S}H^s(\Omega)$. If $s = \infty$, then $k = \infty$. Moreover, $Q$ can be constructed so that $K$ depends continuously on $g$ in a small neighborhood of a fixed $g_0 \in C^k(\bar{\Omega})$.*

*Proof:* Here we will work with $\Psi$DOs with symbols and amplitudes of finite smoothness $k$ with respect to $x$, $y$, and $\xi$, see e.g., [SU2], that satisfy a finite number of the seminorm estimates. A $\Psi$DO of order 0 with such a symbol is bounded in $L^2$ for $k = 2n + 1$, composition of two $\Psi$DOs is a $\Psi$DO with similar symbol for $k \gg 1$, construction of a parametrix to an elliptic $\Psi$DO up to a smoothing operator of finite order also requires finitely many steps and derivatives. Similarly, for any $m$, $s$, a $\Psi$DO of order $m$ is bounded locally as an operator mapping $H^{s+m}$ into $H^s$, provided that its symbol satisfies a finite number of seminorm estimates.

We start with a parametrix of the elliptic operator $W$, see (3.1). In contrast to section 3, it is enough to have the smoothing part $K$ to be of finite order only. Moreover, we have to make the construction uniform for $g$ in a small neighborhood of a given $g_0$. To this end, we will choose all cut-off functions involved in such a construction independent of $g$. Observe that it is enough to work with classical $\Psi$DOs only that have (finite) asymptotic expansions in homogeneous symbols in $\xi$. For $m > 0$, we construct $B = b(x, D)$ to be such that $b \circ \sigma(\Delta^s) = \operatorname{Id} \mod S^{-m}$ near $\Omega_1$, and $b$ is constructed with finitely many iterations, as mentioned above. The symbol $\circ$ stands for composition of symbols by means of finite sums, with enough terms to justify the estimate on the remainder. For any $m > 0$ one has $k > 0$, such that the above construction is possible for $g \in C^k$. We construct $p_0$ such that

$$p_0 \circ \left( \sigma(N) + |\xi|_e^{-1} \circ \sigma(d) \circ b \circ \sigma(\delta) \right) = \operatorname{Id} \mod S^{-m}.$$

Above, $|\xi|_e$ is the Euclidean norm of $\xi$. The symbol of $\mathcal{S}_{\Omega_2}$ mod $S^{-m}$, can be written as $\Lambda = \mathrm{Id} - \sigma(d) \circ b \circ \sigma(\delta)$. Then
$$p_1 := \Lambda \circ p_0 \circ \Lambda$$
satisfies $p_1 \circ \sigma(N) = \Lambda$ mod $S^{-m}$ near $\bar{\Omega}_1$. Moreover, the finitely many seminorms needed for the $H^{s+1} \to H^s$ boundedness of $P_1 := p_1(x, D)$ in any compact in $\Omega_2$ for any fixed $s$, can be estimated by finitely many seminorms of the symbols involved above, which in turn depends on finitely many derivatives of $g$.

So we get
$$P_1 N f = f - d B \delta f + K_1 f \quad \text{in } \Omega_2, \ \forall f \in L^2(\Omega), \tag{5.5}$$
and $K_1 : L^2(\Omega_1) \to H^s(\Omega_1)$, for any fixed $s$, if $g$ is smooth enough; moreover $K_1$ depends continuously on $g$. As a consequence,
$$P_1 N f = f_{\Omega_1}^s + dw + K_1 f \quad \text{in } \Omega_1, \tag{5.6}$$
where $w = v_{\Omega_1} - B\delta f = (\Delta_{\Omega_1, D}^s)^{-1} \delta f - B\delta f$, and $w|_{\partial \Omega_1} = -B\delta f|_{\partial \Omega_1}$. For any $s > 0$, the map $H^s(\Omega) \ni f \mapsto w|_{\partial \Omega_1} \in H^s(\partial \Omega_1)$ is continuous if $g$ is smooth enough because the kernel of $B$ has any fixed number of continuous derivatives away form the diagonal, if $g$ is smooth enough. Moreover, it is continuous in $g \in C^k$, $k \gg 1$. Since $B$ is a parametrix of $\Delta^s$, for any fixed $s$, the map $L^2(\Omega) \ni f \mapsto \Delta^s w \in H^s(\Omega_1)$ is bounded, if $k \gg 1$. Considering $w$ as a solution of a Dirichlet problem in $\Omega_1$, we get that $dw$ in (5.6) belongs to $H^{s+1}(\Omega_1)$, thus we can write (5.6) as
$$P_1 N f = f_{\Omega_1}^s + K_2 f \quad \text{in } \Omega_1, \tag{5.7}$$
where $K_2$ has the properties of $K_1$.

Next, compare $f_\Omega^s$ and $f_{\Omega_1}^s$. We have $f_\Omega^s = f_{\Omega_1}^s + du$ in $\Omega$, where $u = v_{\Omega_1} - v_\Omega$. The vector field $u$ solves the BVP
$$\Delta^s u = 0 \quad \text{in } \Omega, \quad u|_{\partial \Omega} = v_{\Omega_1}|_{\partial \Omega}. \tag{5.8}$$
We need to express $v_{\Omega_1}|_{\partial \Omega}$ in terms of $Nf$. This can be done as follows. By (5.7), and the fact that $f = 0$ outside $\Omega$, one has $-dv_{\Omega_1} = P_1 N f - K_2 f$ in $\Omega_1 \setminus \Omega$. For $(x, \xi)$ in a one-sided neighborhood of $(x_0, \nu(x_0)) \in \Gamma_+$ in $T(\Omega_1 \setminus \Omega)$, where $\nu(x_0)$ is the outer unit normal to $\partial \Omega$, integrate the above along $\gamma_{x,\xi}$ until this geodesic hits $\partial \Omega_1$, where $v_{\Omega_1} = 0$; denote the corresponding time by $\tau(x, \xi)$. We therefore get
$$[v_{\Omega_1}(x)]_i \xi^i = \int_0^{\tau(x,\xi)} [P_1 N f - K_2 f]_{ij}(\gamma_{x,\xi}(t)) \dot\gamma_{x,\xi}^i(t) \dot\gamma_{x,\xi}^j(t) \, dt.$$

Clearly, for any fixed $x$, a set of $n$ linearly independent $\xi$'s in any neighborhood of $\nu(x_0)$ is enough to determine $v_{\Omega_1}(x)$. We choose this set independent of $x$ in a neighborhood of each $x_0 \in \partial \Omega$, then by compactness argument we choose a finite covering and finite number of such sets. This allows us to construct an operator $P_2$, such that
$$v_{\Omega_1}|_{\partial \Omega} = P_2(P_1 N - K_2) f. \tag{5.9}$$
We proved in [SU3] that
$$\|P_2 P_1 h\|_{H^{1/2}(\partial \Omega)} \leq C \|h\|_{\widetilde{H}^2(\Omega_1)}, \quad \forall h \in \widetilde{H}^2(\Omega_1),$$
and our arguments above also show that $P_2 K_2$ depends continuously on $g \in C^k$, $k \gg 1$.

Let $R : H^{s-\frac12}(\partial\Omega) \to H^s(\Omega)$, be the solution operator $u = Rh$ of the boundary value problem

$$\Delta^s u = 0 \quad \text{in } \Omega, \quad u|_{\partial\Omega} = h. \tag{5.10}$$

Lemma 1 and the remark after it imply that $R$ depends continuously on $g$ in the sense of the proposition. Then (5.8) and (5.9) show that $u = RP_2(P_1N - K_2)f$. This and (5.7) yield

$$\begin{aligned} f^s_\Omega &= f^s_{\Omega_1} + du = (P_1N - K_2)f + dRP_2(P_1N - K_2)f \\ &= (\mathrm{Id} + dRP_2)P_1 Nf + Kf, \end{aligned}$$

where $K$ has the properties required. To complete the proof, apply $\mathcal{S}_\Omega$ to the identity above and set $Q = \mathcal{S}_\Omega(\mathrm{Id} + dRP_2)P_1$. $\square$

*Proof of Theorem 2:* By Proposition 4,
$$QN = \mathcal{S} + K.$$

Since $N\mathcal{S} = N$, and $\mathcal{S}Q = Q$, we get that $K = \mathcal{S}K\mathcal{S}$. Apply $\mathcal{S} + K^*$ to both sides above to get

$$\tilde{Q}N = \mathcal{S} + \tilde{K},$$

where $\tilde{Q} = \mathcal{S}(\mathrm{Id} + K^*)Q$, and $\tilde{K} = K^* + K + K^*K$ is self-adjoint and compact in $L^2(\Omega)$, and maps $L^2(\Omega)$ into $\mathcal{S}L^2(\Omega)$. This implies

$$\tilde{Q}N + \mathcal{P} = \mathrm{Id} + \tilde{K} \quad \text{on } L^2(\Omega).$$

Even though a priori we have $\tilde{Q}N : H^1(\Omega) \to L^2(\Omega)$, the proposition shows that $\tilde{Q}N$ extends to a bounded operator on $L^2(\Omega)$. In what follows, we will indicate the dependence on $g$ by placing the subscript $g$ on $N$, etc. Assume that $I_{g_0}$ is s-injective (in $\Omega$) for some simple $g_0 \in C^k(\bar{\Omega})$ and below, assume that $g$ belongs to a small $C^k$ neighborhood of $g_0$. Then $N_{g_0} : L^2(\Omega) \to L^2(\Omega_1)$ is also s-injective. Indeed, assume that $f \in \mathcal{S}L^2(\Omega)$, and $N_{g_0}f = 0$ in $\Omega_1$. Then $\int_\Omega \bar{f} N_{g_0} f \, dx = \int_{\Gamma_-} |I_{g_0}f|^2 \, d\mu = 0$, therefore $f = 0$.

On $L^2(\Omega)$, $\mathrm{Id} + \tilde{K}_{g_0}$ has a finitely dimensional kernel $\mathcal{F}$ of solenoidal tensors, and let $\{f_1,\ldots,f_k\}$ be a basis in it. We can choose it such that $\{N_{g_0}f_1,\ldots,N_{g_0}f_k\}$ is an orthonormal basis in $N_{g_0}\mathcal{F} \subset L^2(\Omega_1)$ because $N_{g_0} : \mathcal{S}_{g_0}L^2(\Omega) \to L^2(\Omega_1)$ is injective. We define the finite rank operator $Q_0 : L^2(\Omega_1) \to L^2(\Omega)$ by setting $Q_0 h = \sum_j (h, N_{g_0}f_j)f_j$, where the inner product is in $L^2(\Omega_1)$. Set

$$Q_g^\sharp = \tilde{Q}_g + Q_0,$$

Then
$$(Q_g^\sharp N_g + \mathcal{P}_g)f = (\mathrm{Id} + K_g^\sharp)f, \quad \forall f \in L^2(\Omega), \tag{5.11}$$

with $K_g^\sharp = \tilde{K}_g + Q_0 N_g$ compact. We claim that $\mathrm{Id} + K_g^\sharp$ is injective for $g = g_0$. Assume that $(\mathrm{Id} + \tilde{K}_{g_0})f + Q_0 N_{g_0}f = 0$. The first term above is orthogonal to $\mathcal{F}$, the second one belongs to $\mathcal{F}$. Therefore, they both vanish, which implies $f \in \mathcal{F}$ and $Q_0 N_{g_0}f = 0$. The explicit form of $Q_0$ yields $(N_{g_0}f, N_{g_0}f_j) = 0$, $\forall j$, and since $f$ is a linear combination of the $f_j$'s, we get $f = 0$.

Therefore, $K_{g_0}^\sharp$ is a compact operator on $L^2(\Omega)$ with $\mathrm{Id} + K_{g_0}^\sharp$ injective. This implies that $\mathrm{Id} + K_{g_0}^\sharp$ is actually invertible. Then (5.11) yields the estimate in Theorem 2 for $g = g_0$:

$$\|f^s\|_{L^2(\Omega)} \leq C\|N_{g_0}f\|_{\tilde{H}^2(\Omega_1)}, \quad \forall f \in H^1(\Omega).$$

To show that $C$ can be chosen independently of $g$ near $g = g_0$, it is enough to observe that $K_g^\sharp$, considered as an operator in $L^2(\Omega)$, depends continuously on $g$ for $k$ large enough. Indeed this is true for $\widetilde{K}_g$ by Proposition 4, and it is also true for $N_g$ in the same space, see for example the representation (3.4). Therefore, for $g$ close enough to $g_0$ in some $C^k$ topology, $\mathrm{Id} + K_g^\sharp$ remains invertible with a uniform bound on the inverse, i.e., $\|f\| \leq C\|(\mathrm{Id} + K_g^\sharp)f\|$ with some $C > 0$ independent on $g$. Let $f \in H^1(\Omega)$, and substitute $f = f_{\Omega,g}^s$ in (5.11) to get

$$\begin{aligned}
\|f_{\Omega,g}^s\|_{L^2(\Omega)} &\leq C\|Q_g^\sharp N_g f\|_{L^2(\Omega)} \\
&\leq C\left(\|\widetilde{Q}_g N_g f\|_{L^2(\Omega)} + \|Q_0 N_g f\|_{L^2(\Omega)}\right) \\
&\leq C'\left(\|N_g f\|_{\widetilde{H}^2(\Omega_1)} + \|N_g f\|_{L^2(\Omega_1)}\right) \leq C''\|N_g f\|_{\widetilde{H}^2(\Omega_1)}, \qquad \forall f \in H^1(\Omega),
\end{aligned}$$

with $C'' > 0$ also independent of $g$ in a neighborhood of $g_0$. This completes the proof that $\mathcal{G}^k$ is open for some $k$.

To complete the proof of Theorem 2, it is enough to observe that analytic functions are dense in $C^k(\bar{\Omega})$, and if $g$ is close enough to a fixed simple metric $g_0$ in $C^2(\bar{\Omega})$, then $g$ is also a metric and also simple. $\square$

## 6 Generic boundary rigidity; proof of Theorem 3

Since Theorem 3 follows from Theorem 5, we will only indicate here the slight changes in the proof of Theorem 5 in [SU3] that imply Theorem 3. The advantage we have here, compared to [SU3], is that we know that the constant $C$ in (1.2) is locally uniform in $g$.

We start with a proposition that allows us to think of classes of isometric metrics in $C^k$, instead of a single metric.

**Lemma 6** *Let $g, \tilde{g} \in C^k(\bar{\Omega})$, $k \geq 1$, and $\psi : \bar{\Omega} \to \bar{\Omega}$ be a $C^l(\bar{\Omega})$ diffeomorphism fixing $\partial\Omega$ with $2 \leq l \leq k+1$. Assume that $\tilde{g} = \psi^* g$. Then $\psi \in C^{k+1}(\bar{\Omega})$, and $\|\psi\|_{C^{k+1}} \leq C(A)$, where $A$ is an upper bound of $\|g\|_{C^k} + \|\tilde{g}\|_{C^k}$.*

*Proof:* We start with the known formula that relates the Christoffel symbols of $g$ and $\tilde{g}$:

$$\Gamma_{ij}^m = \frac{\partial x^m}{\partial \psi^\alpha}\frac{\partial \psi^\beta}{\partial x^i}\frac{\partial \psi^\gamma}{\partial x^j}\tilde{\Gamma}_{\beta\gamma}^\alpha \circ \psi + \frac{\partial x^m}{\partial \psi^\alpha}\frac{\partial^2 \psi^\alpha}{\partial x^i \partial x^j}.$$

Solve this for $\partial^2 \psi^\alpha/\partial x^i \partial x^j$ to get

$$\frac{\partial^2 \psi^m}{\partial x^i \partial x^j} = \frac{\partial \psi^m}{\partial x^\alpha}\Gamma_{ij}^\alpha - \frac{\partial \psi^\beta}{\partial x^i}\frac{\partial \psi^\gamma}{\partial x^j}\tilde{\Gamma}_{\beta\gamma}^m \circ \psi. \tag{6.1}$$

Formula (6.1) was pointed out to the authors by J. Lee. Now, we have $\Gamma_{ij}^m, \tilde{\Gamma}_{ij}^m \in C^{l-1}$, $\nabla \psi \in C^{l-1}$, therefore $\psi \in C^{l+1}$. Iterating this argument, we get that $\psi \in C^{k+1}$. The estimate in the lemma follows immediately for all derivatives of orders between 2 and $k+1$ with $C(A)$ that may a priori depend on $C(A)$, and on a bound of $\nabla \psi$ as well. On the other hand, one can easily get $\|\nabla \psi\|_{C^0} \leq C(A)$ by using the formula for $\psi^* g$. The $C^0$-norm of $\psi$ is bounded by assumption, and this completes the proof of the lemma. $\square$

Let $g_0 \in \mathcal{G}^k$ with $k$ large enough. Let $g_1$ and $g_2$ be two metrics such that $\rho_{g_1} = \rho_{g_2}$ on $\partial\Omega \times \partial\Omega$, and

$$g_1, g_2 \in \mathcal{B} = \left\{ g \in C^k(\bar{\Omega}); \; \|g - g_0\|_{C^k(\bar{\Omega})} \leq \varepsilon \right\} \tag{6.2}$$

We will show that for $0 < \varepsilon \ll 1$, $g_2$ is isometric to $g_1$.

First, by [LSU], there exists a diffeomorphism $\psi$ fixing the boundary, such that $\psi_* g_2$ and $g_1$ coincide at $\partial\Omega$ together with their derivatives up to any fixed order, if $k \gg 1$. The diffeomorphism $\psi$ can be chosen as identifying boundary normal coordinates related to $g_1$ to those related to $g_2$ near $\partial\Omega$, and extending this in the whole domain. Then $\psi = \mathrm{Id} + O(\varepsilon)$ in $C^{k-2}$, therefore, the modified metric $\psi_* g_2$ also belongs to $\mathcal{B}$ with $k$ replaced by $k - 3$, and for some $\varepsilon_1 > 0$, such that $\varepsilon_1 \to 0$, as $\varepsilon \to 0$. Therefore, we may assume that $\psi_* g_2$ is still in $\mathcal{B}$.

Then we pass to semigeodesic coordinates as in Lemma 5, related to each metric, i.e., we replace $g_1$, $\psi_* g_2$ by their push forwards $\phi_{1*} g_1$, $(\phi_2 \circ \psi)_* g_2$ under new diffeomorphisms fixing the boundary. It is important to note that the new metrics still agree at $\partial\Omega$ at any fixed order, if $k \gg 1$ because $\phi_1 = \phi_2$ in $\Omega_1 \setminus \Omega$, see also [SU3]. As above, we can still assume that the new metrics are in $\mathcal{B}$. This gives us that for $f := \phi_{1*} g_1 - (\phi_2 \circ \psi)_* g_2$ we have

$$f \in C^k(\Omega_1), \quad \mathrm{supp}\, f \subset \bar{\Omega}, \quad f_{in} = 0, \; i = 1, \ldots, n. \tag{6.3}$$

We now use the fact that the linearization of $\rho_{g_1}(x, y)$ for $(x, y) \in (\partial\Omega)^2$ is $\frac{1}{2} I_{g_1} f(x, \xi)$ [Sh1] with $\xi = \exp_x^{-1} y / |\exp_x^{-1} y|$, to get

$$\|N_{g_1} f\|_{L^\infty(\Omega_1)} \leq C \|f\|_{C^1}^2, \tag{6.4}$$

with $C$ uniform, if $k \geq 2$. Let $\varepsilon > 0$ be such that $\mathcal{B} \subset \mathcal{G}^k$, and the constant $C$ in (1.2) is uniform in $\mathcal{B}$. Then using (1.2), (6.4), and interpolation estimates, we get that for any $0 < \mu < 1$,

$$\|f^s\|_{L^2} \leq C \|f\|_{L^2}^{1+\mu}$$

with $C > 0$ uniform in $\mathcal{B}$, if $k = k(\mu) \gg 1$. The final step is to estimate $f$ by $f^s$. There is no such estimate for generals $f$'s, but we have the advantage here that $f$ satisfies (6.3). Now, $f_{ni} = 0$ allows us to prove that $\|f\|_{L^2} \leq C \|f^s\|_{H^2}$, see (7.45) and (7.46) in next section. Using interpolation estimates again, we get

$$\|f\|_{L^2} \leq C \|f\|_{L^2}^{1+\mu}$$

with a new $\mu > 0$. This implies $f = 0$, if $\|f\|_{L^2} \ll 1$, and the latter condition is fulfilled, if $\varepsilon \ll 1$. This shows that $g_2 = \phi_* g_1$ with a diffeomorphism $\phi$ that a priori may have lower regularity that $C^{k+1}$. Lemma 6 shows that in fact, $\phi \in C^{k+1}$.

This concludes the sketch of the proof of Theorem 3. $\square$

*Proof of Theorem 4:* Fix $k$ as in Theorem 3. The identity $g_2 = \psi^* g_1$, with $\psi$ as in (1.3), is an equivalence relation, for $k$ fixed, see Lemma 6. Let us denote that relation by $g_1 \sim g_2$. We define $\mathcal{U}$ as follows: it consists of all pairs $(g_1, g_2)$ of simple $C^k$ metrics such that $g_1 \not\sim g_2$ and $g_1, g_2$ have distinct distance functions on the boundary plus all $(g_1, g_2) \in \mathcal{G}^k(\Omega) \times \mathcal{G}^k(\Omega)$ with $g_1 \sim g_2$. We will show first that $\mathcal{U}$ is open in $C^k \times C^k$. Fix $(g_1, g_2) \in \mathcal{U}$. Assume first, that $g_1 \not\sim g_2$. Since the pair is in $\mathcal{U}$, we have that $\rho_{g_1} \neq \rho_{g_2}$ on $(\partial\Omega)^2$. Since $C^k(\bar{\Omega}) \ni g \mapsto \rho_g^2 \in C((\partial\Omega)^2)$ is a continuous map, it follows that a small perturbation of $(g_1, g_2)$ will keep the distance functions distinct. Now, assume that $g_1 \sim g_2$, i.e., $\psi^* g_2 = g_1$, and let $\tilde{g}_1$, $\tilde{g}_2$ be $\varepsilon$-perturbations of $g_1$, $g_2$ in $C^k$. Then $\psi^* \tilde{g}_2$ is an $O(\varepsilon)$ perturbation of $\psi^* g_2 = g_1$ in $C^k$, see Lemma 6. We can apply Theorem 3 to $g_0 := g_1$, and $\tilde{g}_1$, $\psi_* \tilde{g}_2$ (playing the roles of $g_1$ and $g_2$) to conclude that for $0 < \varepsilon \ll 1$, $\tilde{g}_1$, $\psi^* \tilde{g}_2$ are either isometric metrics, both in $\mathcal{G}^k(\Omega)$; or they have distinct boundary distance functions, i.e., $(\tilde{g}_1, \psi^* \tilde{g}_2) \in \mathcal{U}$. This implies the same for $(\tilde{g}_1, \tilde{g}_2)$. Therefore, $\mathcal{U}$ is open.

To prove that $\mathcal{U}$ is dense, observe that any pair $(g_1, g_2)$ of real analytic simple metrics belongs to $\mathcal{U}$. This is true for $g_1 \not\sim g_2$ by the fact that an analytic simple metric is uniquely determined by its boundary distance function, see [LSU]. It is also true when $g_1 \sim g_2$ by Theorem 1.

To prove the last statement of Theorem 4, fix a simple $C^k$ metric $g$. Let $g_2$ be also a simple $C^k$ metric with distinct boundary distance function (then $g \not\sim g_2$). A small perturbation of $g_2$ will preserve this property, as shown above. To prove the density statement, fix a simple $C^k$ metric $g_3$, and $\varepsilon > 0$. If $g_3$ has distinct boundary distance function than that of $g$, we are done. If it is the same, choose $g_4 \in \mathcal{G}^k$, at a distance from $g_3$ not exceeding $\varepsilon/2$. Again, if $g_4$ has distinct boundary distance function than that of $g$, we are done. Otherwise, if $\varepsilon \ll 1$, some $\varepsilon/2$ perturbation $g_5$ of $g_4$ (actually, any that is not equivalent to $g_4$) would necessarily change the boundary distance function of the latter by Theorem 3. Therefore, there is a metric at distance at most $\varepsilon$ from $g_2$ with different boundary distance function than that of $g$. □

# 7 Stability for the non-linear problem

## 7.1 Stability at the boundary

We start with a theorem about stability at the boundary for the non-linear problem. It was first shown by R. Michel [Mi] that the boundary distance function determines all derivatives of $g$ at $\partial \Omega$ in 2 dimensions. In $n \geq 3$ dimensions this was done in [LSU], while the second author and Wang [UW] gave a constructive procedure.

Given two metrics $g_0$, and $g_1$, in a fixed coordinate system, there is a diffeomorphism $\psi$ near $\partial \Omega$ fixing $\partial \Omega$, and mapping the geodesics for $g_0$ normal to $\partial \Omega$ into geodesics for $g_1$ normal to $\partial \Omega$, by preserving the arc-length. Such a diffeomorphism is defined by means of boundary normal coordinates for each metric and can be extended to a global one. Then $g_0$ and $\psi_* g_1$ have common normal geodesics to $\partial \Omega$, close to $\partial \Omega$, and moreover, if $(x', x^n)$ are boundary normal coordinates near a fixed boundary point for one of those metrics, they are also boundary normal coordinates for the other one.

**Theorem 6** *Let $g_0$ and $g_1$ be two simple metrics in $\Omega$, and $\Gamma \subset\subset \Gamma' \subset \partial \Omega$ be two sufficiently small open subsets of the boundary. Let $\psi$ be as above. Then*

$$\left\| \partial_{x^n}^k (\psi_* g_1 - g_0) \right\|_{C^m(\bar{\Gamma})} \leq C_{k,m} \left\| \rho_{g_1}^2 - \rho_{g_0}^2 \right\|_{C^{m+2k+2}(\overline{\Gamma' \times \Gamma'})},$$

*where $C_{k,m}$ depends only on $\Omega$ and on a upper bound of $g_0$, $g_1$ in $C^{m+2k+5}(\bar{\Omega})$.*

*Proof:* In this proof, we will denote $\psi_{1*} g_1$ by $g_1$, thus the normal coordinates $(x', x^n)$ related to $g_0$ are also normal coordinates for $g_1$. It is enough to prove the theorem for $\Gamma$ a small neighborhood of a fixed $x_0 \in \partial \Omega$, and we are going to use boundary normal coordinates there. For any $x = (x', 0)$ close to $(x'_0, 0)$, set $y = (x' + \varepsilon p', 0)$, where $\frac{1}{2} \leq |p'| \leq 2$, and $\varepsilon \geq 0$ is a small parameter. Set

$$\tilde{\rho}_{g_s}(\varepsilon; x', p') = \rho_{g_s}(x, y), \quad s = 0, 1.$$

With some abuse of notation, we will drop the tilde below. Denote

$$\|f\|_m = \|f\|_{C^m(\bar{\Gamma})}, \quad \|\rho^2_{g_1} - \rho^2_{g_0}\|_m = \sup_{|\alpha|+l\leq m,\, x'\in\Gamma',\, \frac{1}{2}\leq|p|\leq 2} \left|\partial^l_\varepsilon\big|_{\varepsilon=0}\partial^\alpha_{x'}(\rho^2_{g_1} - \rho^2_{g_0})\right|.$$

If $[0, 1] \ni t \to \gamma_s(t)$ is the geodesic in the metric $g_s$ connecting $x$ and $y$, with $t$ a natural parameter, then

$$\rho^2_{g_s}(\varepsilon; x', p') = \int_0^1 g_{s,ij}(\gamma_s)\dot\gamma^i_s\dot\gamma^j_s\, dt,$$

and moreover, $\gamma_s$ minimizes the r.h.s. above. We will replace $s \in \{0, 1\}$ by a continuous parameter $s$ by setting $g_s = (1-s)g_0 + sg_1$, $0 \leq s \leq 1$. We show below that $\partial\Omega$ is strictly convex w.r.t. each $g_s$. Then

$$\rho^2_{g_1}(\varepsilon; x', p') - \rho^2_{g_0}(\varepsilon; x', p') = \int_0^1 \frac{d}{ds}\int_0^1 g_{s,ij}(\gamma_s)\dot\gamma^i_s\dot\gamma^j_s\, dt\, ds = \int_0^1 \int_0^1 f_{\alpha\beta}(\gamma_s)\dot\gamma^\alpha_s\dot\gamma^\beta_s\, dt\, ds, \quad (7.1)$$

where $f = g_1 - g_0$, and the Greek symbols vary from 1 to $n - 1$. The terms coming from differentiating $\gamma_s$ above vanish because of the minimizing property of $\gamma_s$ for each fixed $s$. Notice that $(x', x^n)$ are boundary normal coordinates related to $g_s$ for each $s \in [0, 1]$, too, as a consequence of the fact that $(g_s)_{ij} = \delta_{ij}$. Indeed, one can easily verify that for $\varepsilon_0 \ll 1$, the curve $[0, \varepsilon_0] \ni t \mapsto (0, \ldots, 0, 1)$ is a geodesic w.r.t. $g_s$ with $t$ the arc-length; and it is perpendicular to the boundary $x^n = 0$. Introduce

$$I_s f(\varepsilon; x', p') = \int_0^1 f_{\alpha\beta}(\gamma_s)\dot\gamma^\alpha_s\dot\gamma^\beta_s\, dt. \quad (7.2)$$

This is the geodesic X-ray transform that we studied before, related to $g_s$, and written in different coordinates. Then (7.1) can be written as

$$\rho^2_{g_1}(\varepsilon; x', p') - \rho^2_{g_0}(\varepsilon; x', p') = \int_0^1 I_s f(\varepsilon; x', p')\, ds. \quad (7.3)$$

Our goal next is to see that the Taylor expansion of $I_s f$ at $\varepsilon = 0$ determines all derivatives of $f$, and to use (7.3) to derive the same conclusion about the non-linear problem.

In what follows, we work with $g = g_s$ and we will drop the subscript $s$. By the geodesic equation

$$\ddot\gamma^n + \Gamma^n_{\alpha\beta}(\gamma)\dot\gamma^\alpha\dot\gamma^\beta = 0. \quad (7.4)$$

In our coordinates, $\Gamma^n_{\alpha\beta} = -\frac{1}{2}\partial g_{\alpha\beta}/\partial x^n$, and the second fundamental form on $\partial\Omega$ is given by $\Gamma^n_{\alpha\beta}p^\alpha p^\beta$. By the strong convexity assumption, the latter is a positive quadratic form. In particular, we get that the same is true for each $s \in [0, 1]$, with a uniform constant. Set $t = 0$ in (7.4) to get

$$\ddot\gamma^n(0) = -\Gamma^n_{\alpha\beta}(x', 0)\theta^\alpha\theta^\beta,$$

where $\theta = \dot\gamma(0)$. Therefore (see also [Sh1]),

$$\gamma^n(t) = t\theta^n - \frac{t^2}{2}\Gamma^n_{\alpha\beta}(x', 0)\theta^\alpha\theta^\beta + O(t^3),$$

for $\theta$ with a fixed length. In our case we have $\theta = \dot\gamma(0) = \exp^{-1}_{(x',0)}(x' + \varepsilon p', 0) = \varepsilon p' + O(\varepsilon^2)$, so in particular, $|\theta| = O(\varepsilon)$. Replace $\theta$ by $\theta/|\theta|$, and $t$ by $t|\theta|$ above. Then,

$$\gamma^n(t) = t\theta^n - \frac{t^2}{2}\Gamma^n_{\alpha\beta}(x', 0)\theta^\alpha\theta^\beta + O(\varepsilon^3), \quad (7.5)$$

Since $\gamma^n(1) = 0$, we get $\theta^n = \frac{1}{2}\Gamma^n_{\alpha\beta}\theta^\alpha\theta^\beta + O(\varepsilon^3)$. Plug this into (7.5) to get

$$\gamma^n(t) = \frac{\varepsilon^2}{2}\Gamma^n_{\alpha\beta}(x',0)p^\alpha p^\beta(t-t^2) + O(\varepsilon^3). \tag{7.6}$$

Let now $f = f^{(0)}(x') + x^n f^{(1)}(x') + \ldots$ be the Taylor expansion of $f$ near $x^n = 0$. Plug this into (7.2) to get

$$If = \varepsilon^2 f^{(0)}_{\alpha\beta}(x')p^\alpha p^\beta + O(\varepsilon^3). \tag{7.7}$$

This recovers $f^{(0)}_{\alpha\beta}(x')p^\alpha p^\beta$. It is easy to see that any symmetric tensor $h$ can be recovered by $h_{ij} p^i p^j$ known for some finite number $\{p_l\}$ of $p$'s, and moreover, this can be done in a stable way, i.e., we also have $|h| \leq C \sup_l |h_{ij} p^i_l p^j_l|/|p_l|^2$ with $C$ depending on the set. Thus we get

$$\frac{1}{2}\frac{d^2}{d\varepsilon^2}\bigg|_{\varepsilon=0} I_s f = f^{(0)}_{\alpha\beta}(x')p^\alpha p^\beta,$$

and in particular, the l.h.s. above is independent of $s$. This implies the estimate

$$\|f(x',0)\|_m \leq C\|\rho^2_{g_1} - \rho^2_{g_0}\|_{2+m}. \tag{7.8}$$

To study the higher order terms in (7.7), we will plug the Taylor series of $f$ w.r.t. $x^n$ into (7.2), therefore we need first to study integrals of $(x^n)^k f^{(k)}(x')$ over $\gamma$. Observe first that one can expand $\gamma(t)$ into a finite Taylor series in powers of $\varepsilon$, similarly to the second order expansion for $\gamma^n$ in (7.6). Next, $f^{(k)}(\gamma'(t)) = f^{(k)}(x',0) + O(\varepsilon)$, and the remainder can be expanded using higher order derivatives of $f^{(k)}$. So we get

$$\int_0^1 (\gamma^n)^k f^{(k)}_{\alpha\beta}(\gamma')\dot\gamma^\alpha\dot\gamma^\beta\, dt = \varepsilon^{2k+2}(1+O(\varepsilon))\left(\frac{1}{2}\Gamma^n_{\alpha\beta}(x',0)p^\alpha p^\beta\right)^k f^{(k)}_{\alpha\beta}(x',0)p^\alpha p^\beta \int_0^1 (t-t^2)^k\, dt$$

$$= C_k \varepsilon^{2k+2}\left(\Gamma^n_{\alpha\beta}(x',0)p^\alpha p^\beta\right)^k f^{(k)}_{\alpha\beta}(x',0)p^\alpha p^\beta + \varepsilon^{2k+3}\Psi_k, \tag{7.9}$$

where $C_k > 0$, $\Psi_k = \Psi_{k0} + \Psi_{k1}\varepsilon + \Psi_{k2}\varepsilon^2\ldots$, and

$$\|\Psi_{kj}\|_m \leq C\|f^{(k)}\|_{m+j+1}, \tag{7.10}$$

with $C$ above uniform for any fixed collection of indices, if $g_0$ and $g_1$ are bounded in $C^{m+j+2}$.

Consider the fourth order term in (7.7). By (7.9) and (7.2), it involves $f^{(0)}$ and $f^{(1)}$ only and we therefore get

$$\frac{1}{4!}\frac{d^4}{d\varepsilon^4}\bigg|_{\varepsilon=0} I_s f = C_1 \Gamma^n_{\alpha\beta}(x',0)p^\alpha p^\beta f^{(1)}_{\alpha\beta}(x',0)p^\alpha p^\beta + \Psi_{01}, \quad C_1 > 0.$$

We integrate the above in $s$ from 0 to 1, use the fact that $\int_0^1 \Gamma^n_{\alpha\beta} p^\alpha p^\beta ds > c_0|p|^2$ with some $c_0 > 0$ independent of $s$ by the strong convexity assumption, and use (7.10), (7.8) to estimate $\Psi_{01}$ to get

$$\|f^{(1)}\|_m \leq C\|\rho^2_{g_1} - \rho^2_{g_0}\|_{4+m},$$

and $C$ is uniform if $g_0$ and $g_1$ are bounded in $C^{m+2}$.

For general $k \geq 1$, the formula above generalizes to

$$\frac{1}{(2k+2)!}\frac{d^{2k+2}}{d\varepsilon^{2k+2}}\bigg|_{\varepsilon=0} I_s f = C_k \left(\Gamma^n_{\alpha\beta}(x',0) p^\alpha p^\beta\right)^k f^{(k)}_{\alpha\beta}(x',0) p^\alpha p^\beta + \sum_{2k'+j=2k-1} \Psi_{k'j}, \quad (7.11)$$

$C_k > 0$. Now, we can prove the following estimate by induction

$$\|f^{(k)}\|_m \leq C \|\rho^2_{g_1} - \rho^2_{g_0}\|_{2+2k+m}, \quad (7.12)$$

and this estimate requires $m + 2k + 3$ uniformly bounded derivatives of $g_0$ and $g_1$.

To complete the proof of the theorem, it is enough to notice that for a fixed $p'$, $d/d\varepsilon$ is a certain directional derivative w.r.t. $y$, and we need a finite set of $p'$'s. The coordinate change required to pass to the original coordinates, increases the needed number of derivatives of $g$ by 2, which explains the factor $m + 2k + 5$ in the theorem. □

## 7.2 Interior stability, proof of Theorem 5.

Fix $g_0 \in \mathcal{G}^k$, $k \geq k_0$, and let $g$ and $\tilde{g}$ be two metrics as $g_1$ and $g_2$ in Theorem 5 with some $A > 0$ and $\varepsilon_0 \ll 1$, i.e.,

$$\|g\|_{C^k(\bar{\Omega})} + \|\tilde{g}\|_{C^k(\bar{\Omega})} \leq A, \quad \|g - g_0\|_{C(\bar{\Omega})} + \|\tilde{g} - g_0\|_{C(\bar{\Omega})} \leq \varepsilon_0. \quad (7.13)$$

The first condition above is a typical compactness condition. Using the interpolation estimate [Tri]

$$\|f\|_{C^t(\bar{\Omega})} \leq C \|f\|^{1-\theta}_{C^{t_1}(\bar{\Omega})} \|f\|^{\theta}_{C^{t_2}(\bar{\Omega})}, \quad t = (1-\theta)t_1 + \theta t_2, \quad (7.14)$$

where $0 < \theta < 1$, $t_1 \geq 0$, $t_2 \geq 0$, one gets that $\|g - g_0\|_{C^t(\bar{\Omega})} \leq C(A)\varepsilon_0^{(k-t)/k}$ for each $s \geq 0$, if $k > t$; the same is true for $\tilde{g}$. For our purposes, it is enough to apply (7.14) with $t$, $t_1$ and $t_2$ integers only, then (7.14) easily extends to compact manifolds with or without boundary. Set

$$\delta = \|\rho^2 - \tilde{\rho}^2\|_{C(\partial\Omega \times \partial\Omega)}. \quad (7.15)$$

Here and below, a tilde above an object indicates that it is associated with $\tilde{g}$. Using interpolation estimates again, for any $\mu < 1$, we get

$$\|\rho^2 - \tilde{\rho}^2\|_{C^m(\partial\Omega \times \partial\Omega)} \leq C\delta^\mu \quad (7.16)$$

with $C = C(A, \mu, m)$, as long as $k$ is large enough. Here, as in Theorem 6, we prefer to work with the squares of $\rho$, $\tilde{\rho}$ because they are smooth functions with derivatives of any fixed order bounded by $C(A)$, if $k \gg 1$.

In what follows, we denote by $\mu < 1$ constants arbitrarily close to 1, that may change from step to step. We also denote by $C$ various constants depending only on $\Omega$, $A$, $\mu$, and on the choice of $k$ in (7.13). Our goal is to show that for any such $\mu$, there exists $k \gg 1$, and $\varepsilon_0 > 0$, such that the estimate in Theorem 5 holds. We will often use the notation $k \geq k_1(\mu) \gg 1$ to indicate that the corresponding statement holds for $k$ large enough, depending on $\mu$.

By Theorem 6, one can choose a diffeomorphic copy of $\tilde{g}$, that will be denoted by $\tilde{g}$ again, such that the stability estimate in Theorem 6 holds, i.e.,

$$\|\partial^l_\nu(g - \tilde{g})\|_{C^m(\partial\Omega)} \leq C\delta^\mu, \quad \forall l, m, \forall \mu < 1 \quad (7.17)$$

as long as $k \geq k_1(\mu, m + 2l) \gg 1$, where $\partial_\nu$ is the normal derivative. Estimates (7.13) will be replaced by similar ones as in Section 6. Without loss of generality we may assume that the original estimates (7.13) are still satisfied.

Below, we will modify the starting metrics $g$ and $\tilde{g}$ several times, and each subsequent pair will be denoted by $g_l$, $\tilde{g}_l$, where $l = 1, 2, 3, 4$. The corresponding $\rho$'s will be denoted by $\rho_l$, $\tilde{\rho}_l$.

**Construction of $g_1$ and $\tilde{g}_1$.** First, we modify $\tilde{g}$ near $\partial \Omega$ by replacing it there by $g$ in a small $\delta$-dependent neighborhood. Let $\chi \in C^\infty(\mathbf{R})$, such that $\chi(t) = 1$ for $t < 1$, and $\chi(t) = 0$ for $t \geq 2$. Let $M > 0$ be a large parameter that will be specified later. Set

$$\tilde{g}_1 = \tilde{g} + \chi\big(\delta^{-1/M} \rho(x, \partial \Omega)\big)(g - \tilde{g}), \quad g_1 = g. \tag{7.18}$$

Using Taylor's expansion of $g$ and $\tilde{g}$ up to $O\big((x^n)^M\big)$, where $x^n = \rho(x, \partial \Omega)$, and estimate (7.17), we see that for any $m \geq 0$,

$$\|\tilde{g}_1 - \tilde{g}\|_{C^m(\bar{\Omega})} \leq C\delta^{\mu - m/M}, \quad \forall \mu < 1, \tag{7.19}$$

provided that $k \geq k(M, m, \mu)$. We extend $g_1$ and $\tilde{g}_1$ in a small neighborhood $\Omega_1$ of $\Omega$, such that the extended metrics are still simple there, and equal. If $\delta \ll 1$, then (7.13) hold with $A$ and $\varepsilon$ there multiplied by a constant. The modified metrics then satisfy

$$g_1 = \tilde{g}_1 \quad \text{for } -1/C \leq x^n \leq \delta^{1/M}. \tag{7.20}$$

where $x^n$ is the normal coordinate in a collar neighborhood of $\partial \Omega$. In view of (7.19), it is enough to estimate $g_1 - \tilde{g}_1$.

We will use the following observation in what follows. If $g_1$, $\tilde{g}_1$ are $\epsilon$–close in $C^k$ with some $k$, then the corresponding Hamiltonian flows $\Phi^t(x, \xi)$ and $\tilde{\Phi}^t(x, \xi)$ are $O(\epsilon)$ close on any compact set in the $C^{k-2}$ topology w.r.t. the variables $t, x, \xi$. This follows from the fact the Hamiltonian vector fields are $O(\epsilon)$ close in $C^{k-1}$, and if they have $k - 1$ continuous derivatives w.r.t. a parameter, the same is true for the solution (see e.g. [A]) with upper bounds depending on those of the derivatives of order $\geq 1$ of the Hamiltonian field. Now, one can define $g_s$ as above, and choose $s$ to be that parameter, and to apply the mean value theorem for $0 \leq s \leq 1$.

The argument above shows, in particular, that the distance functions $\tilde{\rho}_1$ and $\tilde{\rho}$ related to $\tilde{g}_1$ and $\tilde{g}$, respectively satisfy

$$\|\tilde{\rho}_1^2 - \tilde{\rho}^2\|_{C^{m-2}(\partial \Omega \times \partial \Omega)} \leq C\delta^{\mu - m/M}, \quad \forall \mu < 1, \tag{7.21}$$

for any $m$, if $k \geq k(m, M, \mu)$. To prove this, we write $\rho^2(x, y) = |\exp_x^{-1} y|^2$. Therefore, by choosing $M \gg 1$, we can arrange the estimate (7.16) for $\rho_1^2 - \tilde{\rho}_1^2$ for each fixed $m$ by writing $\rho_1^2 - \tilde{\rho}_1^2 = (\rho^2 - \tilde{\rho}^2) + (\tilde{\rho}^2 - \tilde{\rho}_1^2)$, i.e.,

$$\|\rho_1^2 - \tilde{\rho}_1^2\|_{C^m(\partial \Omega \times \partial \Omega)} \leq C\delta^\mu, \quad \forall \mu < 1, \tag{7.22}$$

as long as $M$ and $k$ are large enough, depending on $m$ and $\mu$.

**Construction of $g_2$ and $\tilde{g}_2$.** Following the proof of Lemma 5, choose $x_0 \in \Omega_1 \setminus \bar{\Omega}$, and let $\psi : \Omega \to W := \psi(\Omega)$ be the corresponding diffeomorphism related to $g_1$. Set also $W_1 = \psi(\Omega_1)$. Denote

$$g_2 = \psi_* g_1, \quad \tilde{g}_2 = \psi_* \tilde{g}_1 \quad \text{in } W. \tag{7.23}$$

Then the straight lines $x' = $ const. are geodesics for $g_2$ but not necessarily for $\tilde{g}_2$. We also have $(g_2)_{in} = \delta_{in}$, $\forall i$.

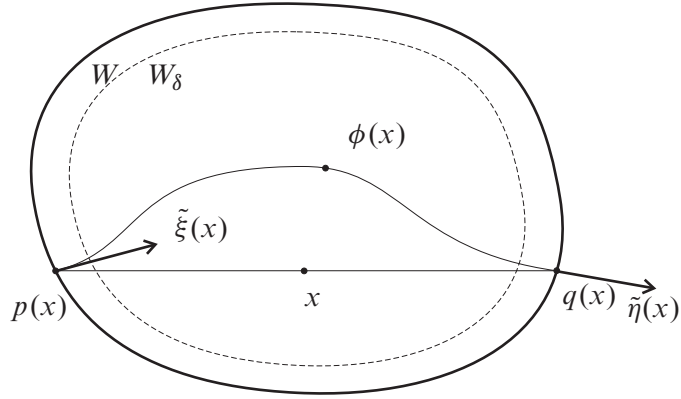

The diffeomorphism $\phi$.

Let $\partial W_{\mp} \subset \partial W$ be the set of those points $x \in \partial W$ with the property that the vector $\pm e_n$ at $x$ points into $W$ (see also (4.5)). For any $x \in W$, let $p(x) \in \partial W_-$, $q(x) \in \partial W_+$ be the endpoints of the line segment in $W$ through $x$ parallel to $e_n$. We will define a new diffeomorphism $\phi$ fixing $\partial W$ in the following way. Set

$$\phi(x) = \widetilde{\exp}_{p(x)} \frac{\tilde{\rho}_2(p(x), q(x))}{|p(x) - q(x)|} |x - p(x)| \tilde{\xi}(x), \quad \tilde{\xi}(x) = -\nabla_p \tilde{\rho}_2(p(x), q(x)), \tag{7.24}$$

where $|\cdot|$ is the Euclidean norm, and $\widetilde{\exp}$ is related to $\tilde{g}_2$. The map $\phi$ is chosen so that $\phi(W) = W$ and $\phi = \text{Id}$ on $\partial W$. Next, $\phi$ also extends into a neighborhood of $W$. Note that $\phi$ maps the straight lines parallel to $e_n$ into geodesics for $\tilde{g}_2$.

The motivation for choosing $\phi$ in such a way comes from an observation made by Michel [Mi]. Since by (7.20), $g_2 = \tilde{g}_2$ on $\partial W$, then $\rho_2 = \tilde{\rho}_2$ on $(\partial W)^2$ would imply $\tilde{\xi}(x) = e_n$. This is implied by the following: by differentiating $\rho_2(x, y) = \tilde{\rho}_2(x, y)$ w.r.t. $x \in \partial W$, for any $y \in \partial W$, we get that the tangential gradients $\nabla' \rho_2$ and $\nabla' \tilde{\rho}_2$ coincide on $(\partial W)^2$. This also allows us to conclude that the full gradients are the same because their lengths are equal (to 1), more precisely, in boundary normal coordinates,

$$\frac{\partial \rho_2}{\partial x^n} = \sqrt{1 - g_2^{\alpha\beta} \frac{\partial \rho_2}{\partial x^\alpha} \frac{\partial \rho_2}{\partial x^\beta}} = \frac{\partial \tilde{\rho}_2}{\partial x^n}. \tag{7.25}$$

In our case, this implies $\tilde{\xi} = e_n$ (under the assumptions $g_2 = \tilde{g}_2$ on $\partial W$, $\rho_2 = \tilde{\rho}_2$ on $(\partial W)^2$). Then $\phi = \phi_0$, where

$$\phi_0 = \widetilde{\exp}_{p(x)} |x - p(x)| e_n = \psi \circ \tilde{\psi}^{-1}. \tag{7.26}$$

Above, $\tilde{\psi}$ is defined as $\psi$ in Lemma 5 but related to $\tilde{g}_1$. The second identity above shows that $\phi_0$ is a diffeomorphism between $W$ and $\phi_0(W)$ that can be extended near $W$ as well. Below we estimate $\phi - \phi_0$ in terms of $\delta$, if $\rho_2$ and $\tilde{\rho}_2$ are not equal on the boundary but satisfy (7.22), also true for $\rho_2, \tilde{\rho}_2$.

We are going to use here (7.20) in a very essential way. For $\delta \ll 1$, set

$$W_\delta = \{x \in W; \ \rho_2(x, \partial W) > \delta^{1/M}/C\}. \tag{7.27}$$

We choose $C > 0$ so that $g_2 = \tilde{g}_2$ for $x \notin W_\delta$, see (7.20). Because of (7.20), the possible singularities in (7.24) connected to $x$, where the ray $s \to x + se_n$ is tangent to $\partial W$ are non-existent because $\phi = \mathrm{Id}$ near such rays.

Observe first that here exists $C > 0$ such that if $|x - y| \leq \delta^{1/2M}/C$ and $x, y \in \partial W$, then $\tilde{\rho}_2(x, y) = \rho_2(x, y)$. We claim that this shows that for any $m$,

$$\left|\partial_x^m\bigl(\tilde{\xi}(x) - e_n\bigr)\right| \leq C\delta^\mu, \quad \left|\partial_x^m\left(\frac{\tilde{\rho}_2(p(x), q(x))}{|p(x) - q(x)|} - 1\right)\right| \leq C\delta^\mu, \quad \forall \mu < 1, \tag{7.28}$$

as long as $M \geq M_1(\mu, m) \gg 1$ and $k \geq k_1(\mu, M, m) \gg 1$. For $m = 0$, the second inequality follows directly from the observation above and (7.22), (7.23) by writing $\tilde{\rho}_2 - \rho_2 = (\tilde{\rho}_2^2 - \rho_2^2)/(\tilde{\rho}_2 + \rho_2)$. To prove the first one for $m = 0$, we need to estimate $\nabla_p\bigl(\tilde{\rho}_2(p(x), q(x)) - \rho_2(p(x), q(x))\bigr)$ for $|p(x) - q(x)| \geq \delta^{1/2M}/C$. We do this for the tangential gradient first using (7.22) again, choosing $M \gg 1$. Then we estimate the remaining normal component of the gradient using (7.25) (note that $x^n$ is a normal coordinate in (7.25) and is not the same as $x^n$ in the fixed coordinate system in $W$ that we are using). This is done by using the estimate for the tangential gradient of $\tilde{\rho}_2 - \rho_2$ and the fact that $|\partial \tilde{\rho}_2/\partial x^n| \geq \delta^{1/2M}/C$, $\partial \rho_2/\partial x^n| \geq \delta^{1/2M}/C$ (which follows by differentiating (7.6)) for $|p(x) - q(x)| < 1/C$ with $C \gg 1$, and by the simplicity condition for $|p(x) - q(x)| \geq 1/C$. This proves (7.28) for $m = 0$. If $|m| > 0$, then (7.28) follows from the arguments above and the fact that subsequent each derivative of $|p(x) - q(x)|^{-1}$, or $(\tilde{\rho}_2(p(x), q(x)))^{-1}$, contributes a factor of the kind $\delta^{-1/2M}$, and a finite product of such factor can be estimated $\delta^{-\kappa}$ for each $\kappa > 0$ small enough, if $M \gg 1$, which can be absorbed by the term $\delta^\mu$.

Those estimates hold near $W$ as well, thus $|\phi - \phi_0| \leq C\delta^\mu$ in a neighborhood of $W$, therefore, we proved that

$$\|\phi - \phi_0\|_{C^m(\bar{W}')} \leq C\delta^\mu, \quad \forall \mu < 1 \tag{7.29}$$

for any $m$, as long as $k \geq k(m, \mu)$ is large enough. In particular, (7.29) shows that $\phi$ is a diffeomorphism for $\delta \ll 1$.

**Construction of $g_3$ and $\tilde{g}_3$.** Set

$$g_3 = g_2, \quad \tilde{g}_3 = \phi^* \tilde{g}_2. \tag{7.30}$$

Then, as mentioned above, the straight lines parallel to $e_n$ are geodesics for $\tilde{g}_3$, with $s$ proportional to the geodesic arc-length but the coefficient of proportionality depends the line. Next, by Lemma 5 and (7.29),

$$(g_3)_{ij} = \delta_{ij}, \quad \|(\tilde{g}_3)_{ij} - \delta_{ij}\|_{C^m(W)} \leq C\delta^\mu, \quad \forall \mu < 1, \tag{7.31}$$

for any fixed $m > 0$, $\mu < 1$, as long as $M \gg 1$, $k \gg 1$.

The new pair of metrics $g_3$ and $\tilde{g}_3$ may not satisfy (7.20) anymore but we will show that they are close in $W \setminus W_\delta$. More precisely, we claim that

$$\sup_{|\alpha| \leq m} \left|\partial^\alpha(\phi(x) - x)\right| \leq C\delta^\mu, \quad \text{for } x \in \bar{W} \setminus W_\delta \tag{7.32}$$

for any $\mu < 1, m$, as long as $M \gg 1$ and $k \gg 1$, depending on $\mu$ and $m$.

To prove (7.32), observe that (7.24) can be written also as

$$\phi(x) = \widetilde{\exp}_{q(x)} \frac{\tilde{\rho}_2(p(x), q(x))}{|p(x) - q(x)|} |q(x) - x| (-\tilde{\eta}(x)), \quad \tilde{\eta} = \nabla_q \tilde{\rho}_2(p(x), q(x)). \tag{7.33}$$

For $\delta \ll 1$, $W_\delta$ is also convex w.r.t. $g_2$. Therefore, for any $x \in W$, we have that the ray $s \mapsto x + se_n$, intersects $\partial W_\delta$ in $m = 0, 1$ or 2 points. If $m = 0$, then along that ray, $g_2 = \tilde{g}_2$ by (7.20), therefore, $\phi(x) = x$. If $m = 1$ or $m = 2$, then either the line segment $[p(x), x]$ is entirely in $W \setminus W_\delta$, or this is true for $[x, q(x)]$. In the first case we use (7.24), in the second one — (7.33). Assume that we have the first case. Then $g_2 = \tilde{g}_2$ near the ray $[x, p(x)]$, therefore $\phi = \phi_0$ near $x$, and (7.32) follows from (7.29). Assume next that we have the second case above. Then we use (7.33), and estimates (7.28) for $\tilde{\xi}$ replaced by $\tilde{\eta}$, and argue as in the proof of (7.29). Note that our choice of $x$ allows us to replace $\widetilde{\exp}$ by $\exp$ in (7.33) in this case.

Estimate (7.32), combined with (7.30), (7.27) implies
$$\sup_{|\alpha| \leq m} |\partial^\alpha (g_3 - \tilde{g}_3)| \leq C\delta^\mu, \quad \text{for } x \in W \setminus W_\delta \tag{7.34}$$

for large $k$ and $M$.

**Construction of $g_4$ and $\tilde{g}_4$.** We will repeat the argument (7.18). Set
$$\tilde{g}_4 = \tilde{g}_3 + \chi(C\delta^{-1/M}\rho_3(x, \partial W))(g_3 - \tilde{g}_3), \quad g_4 = g_3, \tag{7.35}$$

where $C > 0$ is as in (7.27). Then by (7.34), for any $m \geq 0$, $\mu < 1$,
$$\sup_{|\alpha| \leq m} |\partial^\alpha (\tilde{g}_4 - \tilde{g}_3)| \leq C\delta^\mu, \quad \forall x \in W, \tag{7.36}$$

as long as $M \geq M_1(m, \mu) \gg 1$, $k \geq k_1(m, \mu, M) \gg 1$. The advantage that we have with the new metrics $g_4, \tilde{g}_4$, is that they coincide near $\partial W$, i.e.,
$$g_4 = \tilde{g}_4 \quad \text{for } 0 \leq \rho_4(x, \partial W) \leq \delta^{1/M}/C. \tag{7.37}$$

Note that $\rho_4 = \rho_3$ but $\tilde{\rho}_4$ and $\tilde{\rho}_3$ may not be equal, instead they satisfy an estimate similar to (7.21). Therefore,
$$\|\tilde{\rho}_4^2 - \rho_4^2\|_{C^m(\bar{\Omega})} \leq C\delta^\mu, \tag{7.38}$$

for any $m > 0$ and $\mu < 1$ with $M \gg 1$, $k \gg 1$.

**Proof of the stability estimate.** We are ready to linearize the problem now and finish the proof of the stability estimate. Set
$$f = \tilde{g}_4 - g_4. \tag{7.39}$$

Then $f = 0$ near $\partial W$ by (7.37) and we extend it as zero outside of $W$. Arguing as in [SU3], write
$$\tilde{\rho}_4 - \rho_4 = \frac{1}{2} I_{g_4} f + R(f)(x, y), \quad (x, y) \in (\partial W)^2, \tag{7.40}$$

where $y = \exp_x \xi / |\exp_x \xi|$ (the norm here is the same for both metrics because of (7.37)), and the exponential map is in the metric $g_4$. The remainder term satisfies [E]
$$|R(f)(x, y)| \leq C|x - y| \|f\|^2_{C^1(\bar{W})}, \tag{7.41}$$

with $C = C(A)$, and $k = 2$ suffices for this estimate. By (7.37), (7.38), $\tilde{\rho}_4 - \rho_4 = O(\delta^\mu)$, $\forall \mu < 1$, if $M \gg 1$, $k \gg 1$ depending on $\mu$.

We want to apply $I^*_{g_4}$ to both sides of (7.40). One can show that for any simple metric $g$, we have

$$[I^*u]_{ij}(y) = \int_{|\omega|=1} u\left(\gamma_{y,\omega}(\tau_-(y,\omega)), \dot\gamma_{y,\omega}(\tau_-(y,\omega))\right) \omega_i \omega_j \, dS_\omega, \tag{7.42}$$

where $dS_\omega$ is as in Section 2, and $\tau_-(y,\omega) < 0$ is determined by $\gamma_{y,\omega}(\tau_-(y,\omega)) \in \partial W$. Identity (7.42) follows immediately from $(I^*u, h) = (u, Ih)$ and an application of the so-called Santalo formula [Sh1] or arguing as in the proof of Proposition 1 in [SU3].

By (7.42), we have the estimate $\|I^*_{g_4} u\|_{L^\infty(W_1)} \leq C\|u\|_{L^\infty(\Gamma_-)}$. Therefore, (7.40), (7.41) yield

$$\|N_{g_4} f\|_{L^\infty(W_1)} \leq C\delta^\mu + C\|f\|^2_{C^1(\overline W)}, \tag{7.43}$$

for any $\mu < 1$, as long as $k \geq k(\mu)$.

At this point, we apply the stability estimate for the linear problem. By our assumptions, $g_0 \in \mathcal{G}^k$ for $k \gg 1$, and if $\varepsilon_0$ in (7.13) is small enough, then $I_g$ is s-injective in $\Omega$ as well. Then $I_{g_4}$ is s-injective in $W$. By Theorem 2, one has

$$\|f^s\|_{L^2(W)} \leq C\|N_{g_4} f\|_{C^2(W_1)},$$

where we estimated the $\tilde H^2$-norm by the $C^2$ one. The constant $C$ above depends only on $\Omega$, $g_0$, $k$, and $\varepsilon_0$. Recall that $f$ is supported strictly in $W$. Using the interpolation estimate (7.14) again, and continuity properties of the $\Psi$DO $N_{g_3}$, if $k \gg 1$, we get by (7.43),

$$\|f^s\|_{L^2(W)} \leq C\left(\delta^\mu + \|f\|^2_{C^1(\overline W)}\right)^{\mu_1}, \tag{7.44}$$

for any $\mu_1 < 1$, $\mu < 1$, $k \gg 1$, with $k$, $C$ depending on $\mu$, $\mu_1$.

By (7.31), (7.36), $f_{in} = O(\delta^\mu)$, $\forall i$. This estimate and (7.37), (7.39) allow us to estimate $\|f^s\|_{H^2}$ from below by $\|f\|$ modulo $O(\delta^\mu)$. One can express $v$ related to $f$ in terms of $f$ and $f^s$ by solving $dv = f - f^s$, see (4.5). Write $f = f^s + dv$. Then $(dv)_{nn} = \partial_{x^n} v_n = f_{nn} - f^s_{nn}$. Therefore,

$$v_n(x) = \int_{\tau_-(x)}^0 (f_{nn} - f^s_{nn})(x', x^n + s) \, ds, \tag{7.45}$$

where $\tau_-(x) < 0$ is determined by $(x', x^n + \tau_-(x)) \in \partial W$. This yields $\|v_n\|_{L^2(W)} \leq C(\|f^s\|_{L^2(W)} + \delta^\mu)$. We use an interpolation inequality similar to (7.14) but for Sobolev spaces, see [Tri], to get $\|v_n\|_{H^1(W)} \leq C(\|f^s\|_{L^2(W)} + \delta^\mu)^{\mu_2}$, $\forall \mu_2 < 1$, if $k \gg 1$. Next, $\nabla_n v_i + \nabla_i v_n = 2(f_{in} - f^s_{in})$. We write this in the form (4.2), estimate the $L^2$-norm of $v_i$, and then as above, its $H^1$-norm. This yields

$$\|v\|_{H^1(W)} \leq C\left(\|f^s\|_{L^2(W)} + \delta^\mu\right)^{\mu_2}, \quad \forall \mu < 1, \mu_2 < 1.$$

Therefore, for $f = f^s + dv$ we get

$$\|f\|_{L^2(W)} \leq C\left(\|f^s\|^\mu_{L^2(W)} + \delta^\mu\right), \quad \forall \mu < 1. \tag{7.46}$$

We combine this with (7.44). Therefore,

$$\|f\|_{L^2(W)} \leq C\left(\|f\|^{2\mu}_{C^1(\overline W)} + \delta^\mu\right).$$

$\forall \mu < 1$, if $k \gg 1$. We can use interpolation again to replace the $C^1$ norm above by a $C(\bar{W})$ norm; and next, the $L^2$ norm by a $C(\bar{W})$ one. Now, if $1/2 < \mu < 1$, and $\|f\|_{C(\bar{W})} \ll 1$, we get that $\|f\|_{C(\bar{W})} \leq C\delta^\mu$.

We have the following estimates

$$\begin{aligned} f &= \tilde{g}_4 - g_4 = \tilde{g}_3 - g_3 + O(\delta^\mu) & \text{by (7.33), (7.35)} \\ &= \phi^* \tilde{g}_2 - g_2 + O(\delta^\mu) & \text{by (7.30)} \\ &= \phi^* \psi_* \tilde{g}_1 - \psi_* g_1 + O(\delta^\mu) & \text{by (7.23)} \\ &= \phi^* \psi_* \tilde{g} - \psi_* g + O(\delta^\mu). & \text{by (7.18), (7.19)} \end{aligned}$$

By (7.29), (7.26), $\phi^* \psi_* \tilde{g} = \phi_0^* \psi_* \tilde{g} + O(\delta^\mu) = \tilde{\psi}_* \tilde{g} + O(\delta^\mu)$. Hence,

$$f = \tilde{\psi}_* \tilde{g} - \psi_* g + O(\delta^\mu). \tag{7.47}$$

By (7.13), $\tilde{\psi} = \psi + O(\varepsilon_0)$. This and (7.47) show first that $\|f\|_{C(\bar{W})} \ll 1$, if $\varepsilon_0 \ll 1$, see (7.13), therefore, by the arguments above, $\|f\|_{C(\bar{W})} \leq C\delta^\mu$. Then (7.47) again implies $\|\tilde{\psi}_* \tilde{g} - \psi_* g\|_{C(\bar{W})} \leq C\delta^\mu$. The estimate in the $C^2$ norm is obtained by interpolation.

This completes the proof of Theorem 5.